\documentclass{birkjour}
\usepackage{graphicx}%
\usepackage{multirow}%
\usepackage{amsmath,amssymb,amsfonts}%
\usepackage{amsthm}%
\usepackage{mathrsfs,mathtools}%
\usepackage[title]{appendix}%
\usepackage{xcolor}%
\usepackage{textcomp}%
\usepackage{manyfoot}%
\usepackage{booktabs}%
\usepackage{algorithm}%
\usepackage{algorithmicx}%
\usepackage{algpseudocode}%
\usepackage{listings,float}%
 \newtheorem{theorem}{Theorem}[section]
 
 \newtheorem{lemma}[theorem]{Lemma}
 
 \theoremstyle{definition}
 \newtheorem{definition}[theorem]{Definition}
 \theoremstyle{remark}

 \newtheorem{conjecture}{Conjecture}
 \numberwithin{equation}{section}

\begin{document}

%
%
%
%
%
%
%
%
%

\title[Spectral Theory of Hypergraphs]
 {Spectral Theory of Hypergraphs: A Survey}

\author[S. S. Shetty]{Shashwath S Shetty}

\address{%
Department of Mathematics, Manipal Institute of Technology, Manipal Academy of Higher Education, Manipal, Karnataka 576 104. India}

\email{shashwathsshetty01334@gmail.com}

\author[K. A. Bhat]{K Arathi Bhat$^*$\footnote{$^*$Corresponding author}}
\address{Department of Mathematics, Manipal Institute of Technology, Manipal Academy of Higher Education, Manipal, Karnataka 576 104. India}
\email{arathi.bhat@manipal.edu}


\date{\today}

\begin{abstract}
Hypergraphs require higher-dimensional representations, which makes it more difficult to compute and interpret their spectral properties.
This survey article uses the framework of hypermatrices to give an in-depth overview of the spectral theory of hypergraphs.
Our focus in this article relies on the theoretical aspects of hypergraphs that help to ease the computational methods.
Spectral theory hypergraphs, one of the most advanced fields of study, are constantly finding novel applications in various domains such as theoretical computer science, quantum physics, and theoretical chemistry, among many others.
We start our journey by introducing hypergraphs, hypermatrices (tensors), resultants, and their properties. 
We outline some of the results used to determine the adjacency spectrum of hypergraphs and go over some of the groundbreaking findings in the development of the theory.
On passing through a list of bounds for the spectral radius of uniform hypergraphs, we will have a look into the spectral versions of Tur\'an-type problems in hypergraphs.
Finally, in addition to the Estrada index of hypergraphs, some significant results related to the characteristic polynomial and its relation with the matching polynomial are presented.
\end{abstract}

\maketitle

\section{Introduction}
The definitions of the hypermatrix eigenvalues were discussed long back in 1947 \cite{lyusternik1947topological}. 
However, the formal study of the eigenvalues and eigenvectors of the hypermatrices was started in 2005 \cite{qi2005eigenvalues,lim2005singular}. 
Although there were several attempts to define and study the eigenvalues, eigenvectors, and other invariants of the matrices associated with the hypergraph \cite{banerjee2021spectrum, feng1996spectra, friedman1995second, rodriguez2009laplacian, lu2013high, lin2017spectral, sripaisan2022algebraic, vishnupriya2022tensor, portugal2021relating} with varying amounts of success, we can surely pin the article \cite{cooper2012spectra} as the spark for igniting the topic {Spectral Theory of Hypergraphs} via hypermatrices (tensors).

Efforts have been made to study the high-order spectrum (or $p$-spectrum or $h$-adjacency spectrum) of simple graphs way back in 1969 \cite{petersdorf1969spektrum}. The definition was modified in \cite{cvetkovic1980spectra}, and the study has been continued \cite{liu2023high,liu2023erd,liu2024high} by associating the hypermatrices to study the high-order spectrum of simple graphs.
Numerous research studies that addressed the study of the spectra of the Laplacian operator \cite{chung1992laplacian, hu2015laplacian, chan2018spectral}, Laplacian \cite{qi2013h, hu2014eigenvectors, hu2015largest,zheng2024zero,shao2015somespectral}, signless-Laplacian \cite{duan2020largest,duan2020some,xiao2021effect}, generalized adjacency ($A_{\alpha}$) \cite{hou2020alpha,zhou2022some,lin2020alpha,lin2020alpha2,zhang2021spectral} and weighted \cite{wan2022spectra,sun2023spectral,galuppi2023spectral,lin2024abc} hypermatrices associated with the hypergraph can be found in the literature.
The studies on the spectral theory of random hypergraphs \cite{lenz2015eigenvalues,cooper2020adjacency,dumitriu2021spectra,chang2023algebraic}, directed hypergraphs \cite{xie2016spectral,banerjee2017spectrum,liu2018p,li2021alpha,shirdel2023spectral}and oriented hypergraphs \cite{reff2014spectral,duttweiler2019spectra,yu2019signed,wang2022spectral} are also initiated parallelly, albeit they are growing slowly. 
Furthermore, an attempt has been made by employing the hyperdeterminant to generalize the Graham-Pollak tree theorem \cite{tauscheck2024generalizations,cooper2024generalization1,cooper2024generalization2,cooper2024note} from the classical distance to the Steiner distance \cite{hakimi1971steiner} (and related studies can be found in \cite{zheng2025hyperdeterminants,du2025spectral}).
But in this short survey, our focus is to take a deep look into the spectra of a uniform hypergraph using the normalized adjacency hypermatrix that has been defined by Cooper and Dutle \cite{cooper2012spectra}, though the study on the adjacency hypermatrix (not a normalized one) of the uniform hypergraph was initiated a bit earlier \cite{rota2009new}. 
Throughout this article, by adjacency hypermatrix of the hypergraph, we mean the normalized adjacency hypermatrix of the uniform hypergraph \cite{cooper2012spectra}.
We genuinely apologize to the reader for any missed articles, even though we have done our best to gather the pertinent ones.
The article's primary topic will be covered without further ado. 
We use algebra more often than English to define and describe things because it is our belief that it is the best abstract language for mathematicians.


\subsection{Hypergraph}

A hypergraph $\mathcal{H}$ is an ordered pair $(\mathcal{V}:=(\mathcal{V}(\mathcal{H})), \mathcal{E}:=\mathcal{E}(\mathcal{H})),$ where the elements of the set $\mathcal{V}$ are called the vertices and the (possibly multi) set $\mathcal{E} \subseteq \mathcal{P}^*(\mathcal{V})$ contains the non-empty subsets of the vertex set $\mathcal{V}.$ 
A simple hypergraph is a hypergraph in which for any $e,e' \in \mathcal{E},$ neither $e \subseteq e'$ nor $e' \subseteq e.$ 
The \emph{order} of a simple hypergraph $\mathcal{H}$ is $\vert \mathcal{V} \vert$, and the \emph{size} of $\mathcal{H}$ is $\vert \mathcal{E} \vert.$
For a hypergraph, $\max\{ \vert e\vert \vert e \in \mathcal{E} \} $ defines the rank, whereas $\min\{ \vert e\vert \vert e \in \mathcal{E} \} $ defines the co-rank. 
A hypergraph is said to be uniform if $rank(\mathcal{H})=co-rank(\mathcal{H}).$
In other words, a $h$-uniform ($h \geq 2$) hypergraph $\mathcal{H}$ is a hypergraph in which $\vert e \vert=h,$ for each hyperedge $e$ of $\mathcal{H}.$ 
By definition, a $h$-uniform hypergraph is simple if $\mathcal{E}$ is a set (not multi-set). In short, we denote an $h$-uniform hypergraph of order $n$ and size $m$ by $\prescript{n}{m}{H}^{(h)}$. 
For a vertex $u$ in a hypergraph $\mathcal{H}$, the count of hyperedges that contain $u$ is the degree of $u$, and is denoted by $d_{\mathcal{H}}(u)$ or simply $d_u$ (if the hypergraph under consideration is clear). 
The co-degree of two adjacent vertices $u_1$ and $u_2$, denoted by $d_{u_1u_2}$, is the number of hyperedges that contain both $u_1$ and $u_2$. The co-degree of two non-adjacent vertices in hypergraph is equal to zero.
Given two hypergraphs $\mathcal{H}_1=(\mathcal{V}_1,\mathcal{E}_1)$ and $\mathcal{H}_2=(\mathcal{V}_2,\mathcal{E}_2)$, $\mathcal{H}_1$ is called a sub-hypergraph of $\mathcal{H}_2$ if $\mathcal{V}_1 \subseteq \mathcal{V}_2$ and $\mathcal{E}_1 \subseteq \mathcal{E}_2.$
A sub-hypergraph $\mathcal{H}_1$ of $\mathcal{H}_2$ is said to be an induced sub-hypergraph if the edge set $\mathcal{E}_1 = \{ e \in \mathcal{E}_2 \vert e \subseteq \mathcal{V}_1 \}.$ 
Two hypergraphs  $\mathcal{H}=(\mathcal{V},\mathcal{E})$ and $\mathcal{H}'=(\mathcal{V}',\mathcal{E}')$ are claimed to be isomorphic, that is $\mathcal{H} \cong \mathcal{H}'$ if there exist a bijection $\theta$ from $\mathcal{V}$ to $\mathcal{V}'$ such that, a subset $e$ of $\mathcal{V}$ is a hyperedge of $\mathcal{H}$ if and only if $\theta(e) \subseteq \mathcal{V}'$ is a hyperedge of $\mathcal{H}'$.

A walk of length $t \geq 2$ in $\mathcal{H}$ is an alternating sequence $u_0e_1u_1e_2\ldots e_{t} v_t$ of the elements in $\mathcal{V}$ and the elements in $\mathcal{E}$ such that for all $1 \leq i \leq t,$ $u_{i-1},u_i \in e_i$ and $u_{i-1} \neq u_i$. 
If all the vertices and hyperedges in a walk are distinct, then it is called a path.
If there always exists at least one path between every two vertices in $\mathcal{H}$, then $\mathcal{H}$ is connected.
A cycle in $\mathcal{H}$ is a path, but the initial and final vertex coincide.
It can be observed that a simple hypergraph with rank greater than $2$ may contain a cycle of length $2$ (i.e., a $2$-cycle).
A hypergraph is said to be linear if $\vert e \cap e' \vert \leq 1,$ for any two hyperedges $e $ and $e'$ of $\mathcal{H}$. 
A linear hypergraph will not contain a cycle of length two. A hypergraph that is connected and contains no cycle is called a hypertree.
By definition, a hypertree is always linear. The set of all hypertrees of order $n$ (resp. size $m$) that are $h$-uniform is denoted by $\prescript{n}{}{T}^{(h)}$ (resp.  $\prescript{}{m}{T}^{(h)}$).

The number of hyperedges in a path defines its length.
A path of shortest length between two vertices $u$ and $u'$ of a connected component in a hypergraph $\mathcal{H}$ is called a shortest path from $u$ to $u'$.
The distance between $u$ and $u'$ is the length of the shortest path and we denote it by $d_{\mathcal{H}}(u,u')$ or simply $d(u,u')$.
For a vertex $u$ in a connected hypergraph $\mathcal{H}$, the eccentricity is defined as $e_{\mathcal{H}}(u)= \max \{  d(u,v) \vert v \in \mathcal{V}(\mathcal{H})  \}$ and denoted by $e_{\mathcal{H}}(u)$.
The diameter of  $\mathcal{H}$ is defined as $diam(\mathcal{H})=\max\{ e_{\mathcal{H}}(u) \vert u \in \mathcal{V} \}$, and the radius as $rad(\mathcal{H})=\max\{ e_{\mathcal{H}}(u) \vert u \in \mathcal{V} \}$.
A vertex $u$ with $d(u)=1$ is called a pendent or core vertex. In a hypergraph $\mathcal{H}$, if there always exists at least one vertex in every hyperedge of $\mathcal{H}$  that is core, then $\mathcal{H}$ is called a cored hypergraph.

For $i \in [h-1]$, let $V_{e_i} = \{ u_1^{(e_i)}, \ldots, u_h^{(e_i)} \}$ be $h-1$ independent sets of size $h$. 
Then a squid $\mathcal{SQ}^{(h)}$ is a $h$-uniform hypergraph on the vertex set  $ \mathcal{V} = \{ u\} \cup \left(   \cup_{i=1}^{h-1} V_{e_h} \right)$ and the edge set $\mathcal{E} = \{ u, u_1^{(e_1)}, u_1^{(e_2)}, \ldots, u_1^{(e_{h-1})} \} \cup \{ V_{e_i} : i \in [h-1] \}.$ 
 A sunflower hypergraph $\mathcal{S}=\mathcal{S}(m,c,h), ~1 \leq c < h$ is a $h$-uniform connected hypergraph having $m$ hyperedges and each hyperedge of $\mathcal{S}$ contains exactly $c$ vertices of degree $m$ and $h-c$ pendent vertices.
The sunflower $\mathcal{S}(m,1,h)$ is called a hyperstar and is a special class of hypertrees. 
For a non-empty set $P$ and an integer $h \in \mathbb{Z}_+$, we denote the set of all $h$ element subsets of $P$ by $\binom{P}{h}.$ 
A complete $h$-uniform hypergraph on the vertex set $\mathcal{V}$ ($\vert \mathcal{V} =n \vert \geq h$ ) is a hypergraph with the edge set $\mathcal{E} = \binom{\mathcal{V}}{h},$ and is denoted by $K_n^{(h)}.$
A complete $h$-uniform hypergraph on $h+1$ vertices is called a simplex, that is $K_{h+1}^{(h)}$.

The $h$-power hypergraph \cite{hu2013cored,hu2015largest} ($h \geq 3$) of a graph $G$, denoted by $G^{h}$ is a $h$-uniform hypergraph obtained from $G$ by including $h-2$ new vertices of degree one (in addition to two end vertices) in each edge. 
In general, a $h$-uniform hypergraph $\mathcal{H}$ is said to be a power hypergraph, if there exist a graph $G$, such that $G^{h}$ is isomorphic to $\mathcal{H}$. 
Shadow (or a $2$-section) of a hypergraph $\mathcal{H}=(\mathcal{V},\mathcal{E})$ is a simple graph $G=(\mathcal{V}',\mathcal{E}')$ with $\mathcal{V}'=\mathcal{V}$ and $\mathcal{E}'=\{ \{u,v \} :  \{u,v \} \subseteq e \text{ for some } e \in \mathcal{E}  \}.$
Let $u,u'$ be two vertices of a $h$-uniform hypergraph $\mathcal{H}$ and $P:=v_0(=u)e_1v_1e_2\ldots e_sv_s(=u')$ be a path of length $s$ connecting $u$ and $u'$.
Then  $\mathcal{P}_s^{(h)}:=(\cup_{i=1}^s e_i, \{e_i \}_{i=1}^s)$ is a $h$-uniform hypergraph with $s$ hyperedges. 
If $\mathcal{P}_s$ is a linear hypergraph, then we call $\mathcal{P}_s $, a loose path from $u_1$ to $u_2$. 
Let $\mathcal{P}_s$ be a loose path from $u_1$ to $u_2$. 
If $d_{\mathcal{H}}(u_1), d_{\mathcal{H}}(u_2) \geq 3$ and $d_{\mathcal{P}_s}(u) = d_{\mathcal{P}_s}(u) $ for any vertex $u$ in $\cup_{i=1}^s e_i \setminus \{u_1,u_2 \},$ then $\mathcal{P}_s$ is called an internal path of $\mathcal{H}$. 
 If $d_{\mathcal{H}}(u_1) \geq 3$ and $d_{\mathcal{P}_s}(u) = d_{\mathcal{P}_s}(u) $ for any vertex $u$ in $\cup_{i=1}^s e_i \setminus \{ u_1 \},$ then $\mathcal{P}_s$ is called a pendent path of $\mathcal{H}$ at $u_1$.

The applications of the theory of hypergraphs is tremendously broad and deep, including multi-layer network analysis \cite{kivela2014multilayer,battiston2020networks}, signal processing \cite{schaub2021signal,aksoy2020hypernetwork}, molecular chemistry \cite{konstantinova2001application,mann2023ai}, complex social networks \cite{zhu2018social,meenakshi2023mathematical,pitas2016graph}, electrical circuits \cite{li2023influence,baymani2024electrophysiological}, mathematical biology \cite{mukhopadhyay2010multi,lugo2021classification}, and many more.
An initial book on the theory of hypergraphs was by Berge \cite{berge1984hypergraphs} in 1984, while a book with a more algebraic approach, including recent developments, was delivered by Bretto \cite{bretto2013hypergraph}.

\subsection{Hypermatrix}

Let $\mathbb{C}^{[h,n]}$ (or $\mathbb{R}^{[h,n]}$) denote the set of all order\footnote{In the terminologies of hypergraph theory \emph{order} means the number of vertices in the hypergraph, and the \emph{order} of a hypermatrix is the number of indices needed to represent each element of multidimensional array uniquely} $h$, dimension $n$ hypermatrices with complex (real) entries. For the sake of simplicity, we denote the set of all order one hypermatrix $\mathbb{C}^{[1,n]}$ (or $\mathbb{R}^{[1,n]}$) by $\mathbb{C}^n$ (or $\mathbb{R}^n$).
Let $\mathcal{M}=(m_{i_1\ldots i_h}) \in \mathbb{C}^{[h,n]}$ and  $\mathbf{y} :=[y_1, \ldots,y_n] \in \mathbb{C}^n.$
A hypermatrix $\mathcal{M}$ is called a symmetric hypermatrix (or a supersymmetric hypermatrix \cite{kofidis2002best})  if $m_{i_1\ldots i_h} = m_{i_{\sigma(1)}\ldots i_{\sigma(h)}},$ for any permutation $\sigma \in \mathcal{S}_h.$  
The diagonal elements of a hypermatrix $\mathcal{M}$ are the entries of the form $m_{i_1\ldots i_h}$ with $i_1=\cdots = i_h$, and $tr(\mathcal{M}) = \sum\limits_{i=1}^n m_{i\ldots i}$. 
The entries other than the diagonal elements are called the off-diagonal entries.
An order $h$, dimension $n$ hypermatrix is called a diagonal hypermatrix if all the off-diagonal entries are zero. 
A diagonal hypermatrix $\mathcal{M}$ is said to be an identity hypermatrix, denoted by $\mathcal{I},$ if all the diagonal entries are equal to one.
Support of a vector $y$ of length $n$ is the set $I \subseteq [n]$ of all indices $i,1 \leq i \leq n$ such that $y_i \neq 0$.
For order $h$, dimension $n$ hypermatrices $\mathcal{M}_1=m^{(1)}_{i_1\ldots i_h}$ and $\mathcal{M}_2=m^{(2)}_{i_1\ldots i_h}$, we say that $\mathcal{M}_1 \leq \mathcal{M}_2$ if $m^{(1)}_{i_1\ldots i_h} \leq m^{(2)}_{i_1\ldots i_h}$, for all $i_j \in[n],j \in [r].$

\begin{definition}\cite{chang2008perron}
    An order $h$, dimension $n$ hypermatrix $\mathcal{M}=(m_{i_1,\ldots,i_h})$ is said to be reducible, if there exists a proper non-empty index subset $I \subset \{1, \ldots, n \}$ such that $m_{1_1, \ldots, i_h}=0$ for all $i_1 \in I$ and $i_2, \ldots, i_h \notin I.$ If $\mathcal{M}$ is not reducible, then it is termed as an irreducible hypermatrix.
\end{definition}

In \cite{wu2025characterization}, the irreducibility of the hypermatrix was characterized  by the connectivity of associated directed hypergraph.
In \cite{cui2015primitive}, the primitivity \cite{chang2011primitivity} of the hypermatrix (see \cite{horn2012matrix} for the matrix case) was characterized by the reachability of the associated directed hypergraph.

\begin{definition}\cite{friedland2013perron}
   For an order $h$, dimension $n$ hypermatrix $\mathcal{M},$ we associate a $h$-partite graph $G_{\mathcal{M}}=(V,E)$ whose vertex set is the disjoint union $V=\bigcup\limits_{j=1}^h V_j$, where $V_j:=\{ 1, \ldots, n\}, 1\leq j \leq h,$ and an element $\{ {i_t, i_s}\}$ of $ V_t \cup V_s, t \neq s$ is in the edge set $E$ of $G_{\mathcal{M}}$ if and only if  $m_{i_1,\ldots,i_h} >0$ ( $\{ i_t,i_s\} \subset \{i_1,\ldots,i_h\} $) for some $h-2$ indices $\{i_1,\ldots,i_h\} \setminus \{ i_t,i_s\}$ (as multi-sets). Then, the hypermatrix is weakly irreducible if $G_{\mathcal{M}}$ is connected.
\end{definition}
\noindent It can be observed that an \emph{irreducible} hypermatrix is always \emph{weakly irreducible}.

The theory of hypermatrices arises naturally in many places \cite{lim2021tensors} and the explicit application of the analysis of which can be found in \cite{mcconnell2014applications}. 
Chang, Qi, and Zhang \cite{chang2013survey} presented a detailed survey on the spectral theory of non-negative hypermatrices in 2013.
A recent book on the spectral theory hypermatrices is by Qi and Luo \cite{qi2017tensor} in 2017.

\subsection{Resultatnt}

Over any algebraically closed field, the existence of the non-trivial solution for the $k$ homogeneous polynomials in $k$ variables can be discussed through the resultants \cite{gelfand1994discriminants,cox2005using}.
If $f_1,f_2,\ldots,f_k$ is a system of $k$ homogeneous polynomials of (total) degrees $h_1,\ldots,h_k$, respectively in $k$ variables, then this system has a non-zero solution if it satisfies $R_{h_1\ldots h_k} \{f_1,\ldots,f_k\}=0$, where $R$ is the polynomial in the coefficients of $f_1,\ldots, f_k$ called the resultant.

It was shown \cite{grenet2010multivariate} that over any field, the computation of the resultant is $\mathbf{NP}$-hard.
As per our knowledge, in general the Resultants can be computed using the following different approaches \cite{morozov2010new}, that are
 Sylvester Matrix \cite{sylvester1840} (For $h=1$ and for any $n \geq 2$, or $n=2$ and for any $h\geq 2$);
      Bezout's formula \cite{bikker1995bezout};
       Koszul complex \cite{chardin1993resultant,anokhina2009resultant} (Given a non-linear map one can construct a Koszul complex whose determinant is exactly equal to the resultant of the original map);
   Poisson's formula (a recursive approach)\cite{cox2005using,gelfand1994discriminants};
  Contour integrals\cite{morozov2012resultants}.

The algebraic approach to determine the eigenvalues and eigenvectors of the hypermatrices is through the multi-polynomial resultants, where the homogeneous polynomials are obtained from the eigenvalue equations of the hypermatrix. 

Let $\mathbb{C}[y_1,\ldots,y_n]$ be the polynomial algebra over the algebraically closed field $\mathbb{C}$. 
In simple terms, $\mathbb{C}[y_1,\ldots,y_n]$ is the collection of all polynomials in $n$ variables. 
Given an $n$-tuple of non-negative integers  $\alpha:=(a_1,\ldots,a_n)$ and a vector $\mathbf{y}=(y_1,\ldots,y_n)$ of indeterminates, we denote the monomial $\prod_{i=1}^n y_i^{a_i}$ by $\mathbf{y}^{\alpha}$. 
Suppose that $F_j \in \mathbb{C}[y_1,\ldots,y_n] $  is a homogeneous polynomial of degree $d_j$, then it can be written as 
$F_j= \sum_{\vert \alpha \vert = d} P_{j,\alpha} \mathbf{y}^{\alpha},$ where $\vert \alpha \vert = \sum_{i=1}^n a_i.$
In the following definition, we give the axioms to satisfy to call an ordering of monomials a monomial order.

\begin{definition}\cite{cox2005using}
    A relation $>$ is called a monomial order on $\mathbb{C}[y_1,\ldots,y_n],$ if it satisfies the following conditions:
    \begin{itemize}
        \item $>$ is a total ordering relation
        \item $>$ is compatible with multiplication in $\mathbb{C}.$ That is, if $\mathbf{y}^{\alpha_1} > \mathbf{y}^{\alpha_2}$ and $\mathbf{y}^{\alpha}$ be any monomial, then $\mathbf{y}^{\alpha+\alpha_1} > \mathbf{y}^{\alpha+\alpha_2}.$  
        \item $>$ is a well ordering. That is, every non-empty collection of monomials has a smallest element with respect to $>$.
    \end{itemize}
\end{definition}

\begin{theorem}\cite{cox2005using}
    Fix degrees $d_1, \ldots, d_n.$ 
    For $j \in [n],$ consider all monomials $\mathbf{y}^\alpha$  of total degree $d_j$ in $y_1, \ldots, y_n.$ 
    Define a variable $v_{j,\alpha}$ for each such monomial. Then there is a unique polynomial $RES \in \mathbb{Z}[\{v_{j,\alpha} \}]$ with the following properties. 
    \begin{itemize}
        \item For $j \in \{ 1, \ldots, n\},$ if $F_j = \sum\limits_{\vert \alpha \vert = d_j} a_{j, \alpha} \mathbf{y}^{\alpha} \in \mathbb{C}[y_1, \ldots, y_n]$ are homogeneous polynomials of degree $d_j$, then the polynomials have a non-trivial common root in $\mathbb{C}^n$ exactly when $Res(F_1, \ldots, F_n)=0$. (If $RES$ denotes the polynomial in the indeterminates $v_{j, \alpha},$ then $Res$ denotes the evaluation of $RES$ at $a_{j,\alpha}$.)
        \item $Res(x_1^{d_1}, \ldots , x_n^{d_n})=1.$
        \item $RES$ is irreducible, even in $\mathbb{C}[\{v_{j,\alpha}\}].$ 
    \end{itemize}
\end{theorem}

Here we state some properties of the resultant that were used in proving many results.

\begin{theorem}\cite{cox2005using}
    Fix degrees $d_1, \ldots, d_n$. Then for $i \in [n],$ $Res$ is homogeneous in the variables $\{ v_{i,\alpha} \}$ with degree $d_1d_2\ldots d_{i-1}d_{i+1} \ldots d_n.$ Furthermore, the total degree of $Res$ is $d_2\ldots d_n + \sum\limits_{i=2}^n d_1d_2\ldots d_{i-1}d_{i+1} \ldots d_n.$
\end{theorem}

\begin{theorem}\cite{cox2005using}
     Let $F_1, \ldots, F_n \in \mathbb{C}[y_1, \ldots, y_n]$ be the homogeneous polynomials of degrees $d_1, \ldots, d_n$ respectively.
     
\noindent [1] If $i < j,$ $$Res(F_1,\ldots,F_j,\ldots,F_i,\ldots, F_n)= 
 (-1)^{d_0\ldots d_n}Res(F_1,\ldots,F_i,\ldots,F_j,\ldots, F_n).$$
 [2] If $F_i=F_i'F_i''$ is a product of homogeneous polynomials of degrees $d_i'$ and $d_i",$ then 
 $$Res(F_1,\ldots,F_i,\ldots, F_n)=Res(F_1,\ldots,F_i',\ldots, F_n)Res(F_1,\ldots,F_i'',\ldots, F_n)$$
\end{theorem}

\begin{theorem}[B{\'e}zout's theorem]\cite{cox2005using}
    If the equations $F_0=\cdots=F_{n}=0$ are of degrees $d_0,\ldots, d_n$ and have finitely many solutions in the $n$-dimensional projective space $\mathbb{P}^n$, then it has exactly $d_0\ldots d_n$ solutions, counted with multiplicity.
\end{theorem}

\begin{theorem}[Poisson product formula]\cite{cox2005using}
Let $F_0,\ldots,F_n$ be homogeneous polynomials of degrees $d_0,\ldots,d_n$, and 
also let $$\overline{F_i}(x_0,\ldots,x_{n-1}):=F_i(x_0,\ldots, x_{n-1},0), 0\leq i \leq n-1;$$ 
$${f_i}(x_0,\ldots,x_{n-1}):=F_i(x_0,\ldots, x_{n-1},1), 0 \leq i \leq n.$$ If $\mathbf{V}$ denotes the affine variety defined by the polynomials $f_0,\ldots,f_{n-1},$ then
$$ Res(F_0,\ldots,F_n) = Res(\overline{F}_0,\ldots,\overline{F}_{n-1})^{d_n} \prod\limits_{p \in \mathbf{V}} (f_n(p))^{m(p)},$$
where $m(p)$ is the multiplicity of the point $p$ in $\mathbf{V}$.
\end{theorem}

\subsection{Generalized Trace}

For the order $h$, dimension $n$ hypermatrices, the Morozov and Shakirov \cite{morozov2011analogue} have given the formula for $det(\mathcal{I}-\mathcal{M})$ using Schur polynomials in the generalized traces.

\begin{definition}
    Let $A$ be an $n \times n$ matrix with variable entries $A_{ij}.$
    Define $t^{th}$ order trace (for some positive integer $t$) of an order $h$, dimension $n$ hypermatrix $\mathcal{M}$ as 
    $$Tr_t(\mathcal{M}) = (h-1)^{n-1} \sum\limits_{p_1+\cdots+p_n = t} \left( \prod\limits_{i=1}^n \frac{\hat{f_i}^{p_i}}{{(t_i(h-1))}!} \right) tr(A^{t(h-1)}), $$
    where $f_i(y_1, \ldots, y_h) = \left( \mathcal{M}\mathbf{y} \right)_i$ and $\hat{f_i}= f_i\left(\frac{\partial}{\partial A_{i1}}, \frac{\partial}{\partial A_{i2}}, \ldots,  \frac{\partial}{\partial A_{in}}\right).$
\end{definition}

\begin{definition}
    Let $P_0=1$ and for $t >0,$ $t^{th}$ Schur polynomial  $P_t \in \mathbb{Z}[x_1, \ldots, x_h]$ is defined as
    $$P_t(x_1, \ldots, x_h) = \sum\limits_{h=1}^t \sum\limits_{\substack{t_1 + \cdots + t_h=t\\ t_j >0, \ \forall j}}\frac{x_{t_1}\ldots x_{t_h}}{h!}. $$
    In other words, $$exp{\left( \sum\limits_{h=1}^{\infty} x_h z^h\right)}= \sum\limits_{h=1}^{\infty} P_t(x_1, \ldots, x_h) z^h.$$
\end{definition}

\begin{theorem}\cite{morozov2011analogue}
    Let $\mathcal{M}$ be an order $h$, dimension $n$ hypermatrix and $\mathcal{I}$ be the identity hypermatrix (same order and dimension). Then,
    $$ det(\mathcal{I} - \mathcal{A}) =  \sum\limits_{t=1}^{\infty} P_t\left(\frac{Tr_1(\mathcal{M})}{1}, \ldots, \frac{Tr_t(\mathcal{M})}{t}\right) =exp\left(\sum\limits_{t=1}^{\infty} - \frac{Tr_t}{t} \right).$$
\end{theorem}

\begin{theorem}\cite{hu2013determinants}
    Suppose $Tr_t(\mathcal{M})$ denote the $t^{th}$ order trace of an order $h$, dimension $n$ hypermatrix $\mathcal{M}$. Then    $Tr_t(\mathcal{M}) = \sum\limits_{\lambda \in \sigma(\mathcal{M})} \lambda^{t},$
    $\sigma(\mathcal{M})$ is the multi-set of the roots of the characteristic polynomial of $\mathcal{M}.$
\end{theorem}

For an integer $t > 0,$ define 
$$ \mathcal{S}_t = \{ ((j_1, \beta_1),\ldots, (j_t, \beta_t)) \vert j_1 \leq \cdots \leq j_t, \beta_i \in  [n]^{h-1}, 1 \leq i \leq d \}. $$

Given an element $S  =((j_1, \beta_1),\ldots, (j_t, \beta_t))\in \mathcal{S}_t$. We call each $(j_i, \beta_i),
1\leq i \leq t,$ a component of $S$. Also, let $j_i, 1 \leq i \leq t$ be the primary element of $S$, and the components of $\beta_i$ be the secondary elements of $S$. 
An element $S \in \mathcal{S}_t$ is said to be $h$-valent, if every $i \in [n]$ occurs exactly $0\mod{h}$ times as some element of $S$.
Also, let $$S_t'= \{S \in S_t \vert S \text{ is } h-\text{valent} \}.$$

\begin{definition}\cite{shao2015some}
    Let $\overrightarrow{G}(M)$ be the weighted digraph associated with an $n \times n$ square matrix $M=(m_{ij})$, that is defined as,
   $ V(\overrightarrow{G}) := [n]$  and there will be an arc $(i,j) \text{ with weight } m_{ij}$ in $\overrightarrow{G}$ if and only if $m_{ij} \neq 0.$
\end{definition}

In other words, if all the non-zero entries of the matrix $M$ are positive integers, the weighted digraph can be instead called a multi-digraph and the weighted arc $(i,j)$ with weight $m_{ij} (\in \mathbb{Z}_+)$ instead called a multi-arc with multiplicity $m_{ij}$.

\begin{definition}\cite{shao2015some}
   Let $S =((j_1, \beta_1),\ldots, (j_t, \beta_t))\in \mathcal{S}_t$, where $j_i \in  [n], \ \beta_i \in [n]^{h-1}$ for $1 \leq i \leq t.$ Then
\begin{itemize}
    \item $E(S) := \bigcup\limits_{k=1}^n E_k(S), $ where $E_k(S):= \{(j_k, u_2), \ldots, (j_k, u_h) \}$ and $\beta_k = u_2\ldots u_h.$
    \item $b(S)$ denotes the product of the factorials of the multiplicities of all the arcs of $E(S)$. 
    \item $c(S)$ denotes the product of the out-degrees of all the vertices that are incident with some arcs of $E(S)$.
    \item $\mathbf{W}(S)$ denotes the set of all eulerian cycles $W$ with arc multi-set $E(W)=E(S).$
\end{itemize}
\end{definition}

\begin{theorem}\cite{shao2015some}
    Let $\mathcal{M}=(m_{i_1\ldots i_h})$ be an order $h$, dimension $n$ hypermatrix, and $Tr_t(\mathcal{M})$ be the $t^{th}$ order trace of $\mathcal{M}.$ Then 
    $$Tr_t(\mathcal{M}) = (h-1)^{n-1} \sum\limits_{S \in \mathcal{S}_t'} \frac{b(S)}{c(S)} \pi_{S}(\mathcal{M}) \vert \mathbf{W}(S) \vert,$$
    where $\pi_{S}(\mathcal{M}) = \prod\limits_{i=1}^n m_{j_i\beta_i}$ and $S=((j_1,\beta_1), \ldots, (j_t,\beta_t)).$
\end{theorem}

In the above result, note that the summation is not taken over all $S \in \mathcal{S}_t $, but only over all $S \in \mathcal{S}_t'.$
This is because\cite{shao2015some} if $ \mathbf{W}(S) \neq \emptyset,$ then $S$ is $h$-valent.



\section{Some results on the eigenvalues}

 Given two hypermatrices $\mathcal{M}_1 =(m_{i_1\ldots i_{h_1}}^{(1)})$ and $\mathcal{M}_2=(m_{i_1\ldots i_{h_12}}^{(2)})$ (both) of dimension $n$, and order $h_1 (h_1 \geq 2)$ and $h_2 (h_2 \geq 1)$, respectively. The (general) product (see \cite{shao2013general})  $\mathcal{M}_1\mathcal{M}_2$ is a dimension $n$ hypermatrix $\mathcal{M}_3$ of order $(h_1-1)(h_2-1)+1$, whose entries are given by 
    $$ m^{(3)}_{i \alpha_1 \ldots \alpha_{h_1-1}}=
    \sum\limits_{i_2,\ldots,i_{h_1} \in [n]} m^{(1)}_{ii_2\ldots i_{h_1}} 
    m_{i_2\alpha_1}^{(2)} \ldots 
    m_{i_{h_1}\alpha_{h_1-1}}^{(2)}, $$ 
    $i \in [n]$ and $\alpha_j \in [n]^{h_2-1}, 2 \leq j \leq h_1.$ 
    A symmetric hypermatrix $\mathcal{M}$ is said to be positive definite if the homogeneous polynomial 
$$\mathcal{M}\mathbf{y}^{h} = \sum\limits_{\substack{i_j \in [n] \\ 1 \leq j \leq h} } m_{i_1 \ldots i_h} y_{i_1}\cdots y_{i_h} $$ 
(of total degree $h$) is positive definite. $\mathcal{M} \mathbf{y}^{h-1}$ is a column vector of dimension $n$ whose $i^{th}$ component is given by  
\begin{equation*}
    (\mathcal{M} \mathbf{y}^{h-1})_i = \sum\limits_{\substack{i_j \in [n] \\ 2 \leq j \leq h} } m_{i i_2 \ldots i_h} y_{i}\cdots y_{i_h}.
\end{equation*}
The symmetric hyperdeterminant \cite{qi2005eigenvalues} of a hypermatrix is the resultant  of the polynomial system $\mathcal{M} \mathbf{y}^{h-1}.$ Later, Hu et al. generalized   the the definition of the determinant, to nonsymmetric hypermatrices (see  \cite{hu2013determinants}).
A complex number $\lambda$ is called an eigenvalue of a hypermatrix $\mathcal{M} \in \mathbb{C}^{[h,n]} $ corresponding to an eigenvector $\mathbf{y} \in \mathcal{C}^n$ if they satisfy, 
$$\mathcal{M}\mathbf{y}^{h-1} = \lambda \mathbf{y}^{[h-1]},$$ 
where $\mathbf{y}^{[h-1]}$ is a column vector of dimension $n$ whose $i^{th}$ component is $y_i^{h-1}$. The multiset of all eigenvalues of $\mathcal{M}$ is called the spectrum of $\mathcal{M},$ denoted by $\sigma(\mathcal{M}).$
A real eigenvalue $\lambda \in \sigma(\mathcal{M})$ is called an $H$-eigenvalue of $\mathcal{M}$ if the eigenvector of $\mathcal{M}$ associated with $\lambda$ is real.
The set of all $H$-eigenvalues of $\mathcal{M}$ are called the $H$-spectrum of $\mathcal{M}$.

If $\mathcal{M}$ is a diagonal hypermatrix of order $2k, k \in \mathbb{N}$, then the diagonal elements of $\mathcal{M}$ are the eigenvalues of $\mathcal{M}$.
The characteristic polynomial of the adjacency hypermatrix of the hypergraph $\mathcal{H}$ is termed the characteristic polynomial of the hypergraph $\mathcal{H},$ and similarly the eigenvalues, spectral radius etc.
In \cite{qi2005eigenvalues}, it has been shown that the eigenvalues of $\mathcal{M}$ are exactly the roots of the polynomial (in $\lambda$) $det(\mathcal{M}-\lambda \mathcal{I}).$ 
Also, the product of the eigenvalues of a hypermatrix $\mathcal{M}$ is equal to $det(\mathcal{M}),$ and the sum of the eigenvalues of $\mathcal{M}$ is equal to $(h-1)^{n-1} tr(\mathcal{M}).$ 
Although, deciding that a given rational number is an eigenvalue (over $\mathbb{R}$) of the hypermatrix is $\mathbf{NP}$-hard and many other hypermatrix problems are $\mathbf{NP}$-hard \cite{hillar2013most}, a lot of computational techniques \cite{cui2014all,chen2015finding, pearson2015spectral,chen2016computing,chen2017computing,nie2018real,chen2018semidefinite, kuo2018continuation,sheng2019local,chen2020trust,ma2023noda,liu2022exact,chen2019homotopy,clark2020stably,chang2023adaptive} have been developed for determining the eigenvalues (also the eigenpair) of a hypermatrix in special cases.

\begin{definition}
 Given two hypergraphs $\mathcal{H}_1=(\mathcal{V}_1, \mathcal{E}_1)$ and $\mathcal{H}_2=(\mathcal{V}_2, \mathcal{E}_2)$, 
    \begin{itemize}
        \item[1] the union of $\mathcal{H}_1$ and $\mathcal{H}_2$, denoted by  $\mathcal{H}_1 \cup \mathcal{H}_2$ is a hypergraph with vertex set $\mathcal{V}_1 \cup \mathcal{V}_2$ and edge set $\mathcal{E}_1 \cup \mathcal{E}_2;$
        \item[2] the Cartesian product of $\mathcal{H}_1$ and $\mathcal{H}_2$, denoted by  $\mathcal{H}_1 \square \mathcal{H}_2$ is a hypergraph with vertex set $\mathcal{V}_1 \times \mathcal{V}_2$ and the edge set is given by 
        $\mathcal{E}(\mathcal{H}_1 \square \mathcal{H}_2)=\{e_i \times \{v_j \}| e_i \in \mathcal{E}_1, v_j \in \mathcal{V}_2\} \cup \{\{ v_i\} \times e_j | v_i \in \mathcal{V}_2, e_j \in \mathcal{E}_2\};$
        \item[2] the tensor product of $\mathcal{H}_1$ and $\mathcal{H}_2$, denoted by  $\mathcal{H}_1 \otimes \mathcal{H}_2$ is a hypergraph with vertex set $\mathcal{V}_1 \times \mathcal{V}_2$ and the edge set is given by 
        $\mathcal{E}(\mathcal{H}_1 \otimes \mathcal{H}_2)=\{ \{ (u_1,v_1),\ldots,(u_h, v_h) \}  : \{ u_1,\ldots,u_h \} \in \mathcal{E}_1\} \text{ and } \{ v_1,\ldots,v_h \} \in \mathcal{E}_2\}.$
    \end{itemize}
\end{definition}

\begin{theorem}\cite{cooper2012spectra}
    Let $\mathcal{H}_1$ and $\mathcal{H}_2$ be  $h$-uniform, vertex disjoint hypergraphs of order $n_1$ and $n_2$ respectively, and $\mathcal{H}=\mathcal{H}_1 \cup \mathcal{H}_2$ (union of two hypergraphs). Then\\
    $$\Phi_{\mathcal{H}}(\lambda)=\Phi_{\mathcal{H}_1}(\lambda)^{(h-1)^{n_2}}\Phi_{\mathcal{H}_2}(\lambda)^{(h-1)^{n_1}}.$$
    In other words, if $\lambda$  is an eigenvalue of $\mathcal{H}_2$ with multiplicity $m,$ then $\lambda$ is an eigenvalue of $\mathcal{H}$ with multiplicity $m^{(h-1)^{n_1}}.$
\end{theorem}

In the above result, the hypergraph $\mathcal{H}$ can be seen as a disconnected hypergraph with components $\mathcal{H}_1$ and $\mathcal{H}_2.$ 
In this direction, the result is generalized to arbitrary symmetric weakly reducible hypermatrix \cite{hu2013determinants,shao2013some}.

\begin{theorem}\cite{cooper2012spectra,pearson2013eigenvalues}
    Let $\lambda_1$ and $\lambda_2$ respectively, be the eigenvalues of two $h$-uniform hypergraphs $\mathcal{H}_1$ and $\mathcal{H}_2$. Then
    \begin{itemize}
        \item $\lambda_1+ \lambda_2$ is an eigenvalue of $\mathcal{H}_1 \square \mathcal{H}_2.$
        \item $\lambda_1\lambda_2$ is an eigenvalue of $\mathcal{H}_1 \otimes \mathcal{H}_2.$
    \end{itemize}
\end{theorem}

The spectrum is said to be $h$-symmetric if it is invariant under the multiplication of any $h^{th}$ root of unity.
A $h$-uniform hypergraph is said to be an $h$-partite hypergraph if there exist an $h$ partition of the vertex set such that each hyperedge contain exactly one vertex from each partition. 
In \cite{cooper2012spectra}, authors showed that the spectrum of the $h$-uniform, $h$-partite hypergraph is $h$-symmetric and also gave an example to show that it is a one-way result, unlike the graph case. Also, they posed a problem of finding more classes of hypergraphs whose spectrum is $h$-symmetric.

In \cite{hu2014eigenvectors}, authors defined a class of bipartite hypergraphs known as $hm$ bipartite hypergraphs. That is, a $h$-uniform bipartite hypergraph with partite set $\mathcal{V}_1$ and $\mathcal{V}_2$ such that every hyperedge of the hypergraph has exactly one vertex from $V_1$ (in other words, exactly $h-1$ from $V_2$). Later, it was generalized \cite{shao2015some} to $p-hm$ bipartite hypergraphs, where every hyperedge of the hypergraph has exactly $p (p > 0)$  vertices from $\mathcal{V}_1.$ 

    \begin{itemize}
        \item\cite{cooper2012spectra} If $\mathcal{H}$ is isomorphic to a $h$-uniform $h$-partite hypergraph, then $\sigma(\mathcal{H})$ is $h$-symmetric.
        \item\cite{hu2014eigenvectors}If $\mathcal{H}$ is isomorphic to a $h$-uniform $hm$ bipartite hypergraph, then $\sigma(\mathcal{H})$ is $h$-symmetric.
        \item\cite{shao2015some} If $\mathcal{H}$ is isomorphic to a $h$-uniform $p-hm$ bipartite hypergraph with $gcd(p,h)=1$, then $\sigma(\mathcal{H})$ is $h$-symmetric.
    \end{itemize}

A hypergraph spectrum is said to be symmetric, if $\lambda \in\sigma(\mathcal{M})$ implies $-\lambda\in \sigma(\mathcal{M})$ with the same multiplicity. 
Pearson and Zhang \cite{pearson2014spectral} proposed the problem of characterizing the $h$-uniform hypergraphs with symmetric spectrum.
Unlike the definition of the symmetric spectrum, the $H$-spectrum is said to be symmetric \cite{nikiforov2017hypergraphs} if $\lambda$ is in $H$-spectrum of $\mathcal{M}$ then $-\lambda$ is also in $H$-spectrum of $\mathcal{M}$ (need not be with same multiplicity). Nikiforov introduced the concept of geo-symmetric spectrum (see \cite{nikiforov2017hypergraphs}) of hypermatrices, where the condition is imposed on the eigenvectors also. Further, he defined the odd-colorable even order hypermatrices and characterized geo-symmetric (so, the symmetric) spectrum using it.
Later, Fan et al. generalized \cite{fan2019spectral} the concept of symmetric spectrum to $t$-symmetric, where the multiplication by a $t^{th}$ root of unity does not affect the multi-set of spectrum. 
Further research in this direction can be found in \cite{hu2020symmetry,fan2021cyclic,fan2022stabilizing,yuan2018some,fan2020minimal,song2024spectral}.

Two hypergraphs $\mathcal{H}_1$ and $\mathcal{H}_2$ are said to be co-spectral if the multi-set of the spectrum $\mathcal{H}_1$ and $\mathcal{H}_2$ are same, and are said to be weakly co-spectral if the underlying set of the spectrum are equal (i.e., their multiplicities need not be same). 
We call a hypergraph $\mathcal{H}$ is $\mathbf{DS}$, if $\sigma(\mathcal{H})$ determines $\mathcal{H}$ uniquely. 
Bu, Zhou and Wei have shown \cite{bu2014cospectral} that the following $h$-uniform hypergraphs are $\mathbf{DS}$:
\begin{itemize}
    \item Complete hypergraph $\mathcal{K}_n^{(h)}$; an edge deleted sub-hypergraph of the complete hypergraph  $\mathcal{K}_n^{(h)} \setminus e$; a $h$-uniform sub-hypergraph of the simplex $\mathcal{K}_{h+1}^{(h)}$; union of simplex $\mathcal{K}_{h+1}^{(h)}$ and $\overline{\mathcal{K}_{n'}}$; $\mathcal{K}_n^{(h)}\setminus \mathcal{E}(\mathcal{S}),$ where $\mathcal{S}$ is a $h$-uniform sub-hypergraph of $\mathcal{H}$ isomorphic to a sunflower $\mathcal{S}(m,h-1,h) $ and $m\leq  n-r.$ 
\end{itemize}

If $\sigma(G_1^{h}) = \sigma(G_2^{h})$, then the simple graphs $G_1$ and $G_2$ are said to be $h$-order $(h \geq 2)$ co-spectral.
Moreover, they are said to be high-order co-spectral if $G_1$ and $G_2$ are $h$-order co-spectral for each $h \geq 2$.
The high-order spectrum of a graph $G$ determines $G$ completely, i.e., $G$ is $\mathbf{DHS}$ if $G_1$ is a graph that is high-order co-spectral with $G$, then $G_1 \cong G$.

Smith graphs are the ones with the spectral radius of at most two. They are called so, as it was completely characterized by Smith \cite{smith1970some} in 1970.
Later, it was shown \cite{van2009developments} that all but two of the Smith graphs are $\mathbf{DS}$.
Schwenk has given a method to construct infinite pairs of co-spectral trees (that are not isomorphic) \cite{schwenk1973almost} and proved that most of the trees are \emph{not} $\mathbf{DS}$. 
Chen, Sun, and Bu \cite{chen2025high} have shown that if $G$ is a Smith graph, then $G$ is $\mathbf{DHS}$, and all the co-spectral trees constructed by Schwenk are $\mathbf{DHS}$. 
He also conjectured that the high-order spectra completely determine all simple graphs.

\begin{theorem}\cite{zhou2014some}
    A non-zero complex number $\lambda \in \sigma(G)$, then $\lambda^{\frac{2}{h}} \in \sigma(G^{h}).$ Moreover, $\rho(G^{h})=\rho(G)^{\frac{2}{h}}.$ 
\end{theorem}

Later, the authors of the article \cite{cardoso2020spectrum} proposed that all the eigenvalues of the power hypergraph (generally, for generalized power) $G^{h}$ can be obtained from the eigenvalues of the subgraph (induced subgraph, if $h=3$) of $G$ by $\beta^2=\lambda^h$. 
In 2023 Chen et al.\cite{chen2023all}, gave a  counter example and proved that it is not just a subgraph, but the signed subgraph. 
A signed graph $G_{\pi}$ is a triple $(V,E,\pi)$ with  $V$ as the vertex set, $E$ as the edge set and a sign function $\pi$ that assigns either $+1$ or $-1$ to each of its edges.

\begin{theorem}\cite{chen2023all}
    If $G^{h}$ denotes the $h$-power of a graph $G$, then a complex number  $\lambda \in \sigma(G^{h})$ if and only if 
    \begin{itemize}
        \item[1.] some \emph{signed induced subgraph} of the graph $G$ has an eigenvalue $\beta$ such that $\beta^2=\lambda^h$, for $h=3;$
        \item[2.] some \emph{signed subgraph} of the graph $G$ has an eigenvalue $\beta$ such that $\beta^2=\lambda^h$, for $h \geq 4.$
    \end{itemize}
\end{theorem}

\subsection{A conjecture on the multiplicities of the eigenvalues}

The multiplicity of an eigenvalue $\lambda_j$ in $\sigma(\mathcal{M})$ is the algebraic multiplicity of $\lambda_j$ in $\mathcal{M}$ and is denoted by $am_{\mathcal{M}}(\lambda_j).$ 
The dimension of the affine variety $V_{\mathcal{M}}(\lambda_j) = \{ \mathbf{y} \in \mathbb{C}^{n} | \mathcal{M}\mathbf{y}^{h-1}= \lambda_j \mathbf{y}^{h-1} \} \subseteq \mathbb{C}^n$ is the geometric multiplicity of $\lambda_j$ in $\mathcal{M}$. 
By $am_{\mathcal{H}}(\lambda)$ (resp. $gm_{\mathcal{H}}(\lambda)$) we mean  $am_{\mathcal{A}(\mathcal{H})}(\lambda)$ (resp. $gm_{\mathcal{A}(\mathcal{H})}(\lambda)$).

\noindent In 2016, Hu and ye \cite{hu2016multiplicities} proposed the following conjecture relating the algebraic and geometric multiplicity of an eigenvalue of a hypermatrix.

\begin{conjecture}\cite{hu2016multiplicities}
    If $\mu$ is an eigenvalue of an order $h$ hypermatrix $\mathcal{M}$ of dimension $n$, then
    \begin{equation*}
    am_{\mathcal{M}}(\mu) \geq gm_{\mathcal{M}}(\mu)(h-1)^{gm_{\mathcal{M}}(\mu)-1}.
    \end{equation*}
\end{conjecture}

\noindent Cooper and Fickes \cite{cooper2021geometric} proved the above conjecture for the zero eigenvalue of a $h$-uniform hyperpath.
\begin{theorem}\cite{cooper2021geometric}
    If  $\mathcal{P}:=\mathcal{P}_m^{(h)}$ denotes a hyperpath with $m$ hyperedges, then 
    \begin{equation*}
        am_{\mathcal{P}}(0) \geq  gm_{\mathcal{P}}(0)(h-1)^{gm_{\mathcal{P}}(0)-1}.
    \end{equation*} 
\end{theorem}
 Zheng \cite{zheng2024algebraic} has extended the above result for any eigenvalue $\lambda$ of a $h$-uniform hyperpath and $h$-uniform hyperstar. 
\begin{theorem}\cite{zheng2024algebraic}
    Suppose that $\mathcal{H}$ is either a hyperpath ($\mathcal{P}_m^{(h)}$) or a hyperstar ($\mathcal{S}_m^{(h)}$). If a complex number $0 \neq \lambda \in \sigma(\mathcal{H})$, then 
    \begin{equation*}
        am_{\mathcal{H}}(\lambda) \geq   gm_{\mathcal{H}}(\lambda) (h-1)^{gm_{\mathcal{H}}(\lambda)-1}.
    \end{equation*} 
\end{theorem}

Fan \cite{fan2024multiplicity} has proved the equality case in the  conjecture for some eigenvalues of the adjacency and Laplacian hypermatrices.
Further developments in proving the above conjecture can be found in \cite{fan2024multiplicity}.

\section{Spectral radius}

    For an order $h$ hypermatrix $\mathcal{M}$ of dimension $n$, the spectral radius is  denoted by $\rho(\mathcal{M})$, and is defined as 
    $$ \rho(\mathcal{M}) = \max \{ \vert \lambda \vert : \mathcal{M}\mathbf{y}^{h-1} = \lambda \mathbf{y}^{[h-1]}, \mathbf{0} \neq \mathbf{y} \in \mathbb{C}^n  \}. $$
For a non-negative symmetric hypermatrix $\mathcal{M}$ it has been shown \cite{qi2005eigenvalues} that 
$$\rho(\mathcal{M}) = \max\limits \{  \mathcal{M} \mathbf{y}^{h} \vert {\mathbf{y} \in \mathbb{R}^n,  ||\mathbf{y}||_h =  1} \}, \text{ where } ||\mathbf{y}||_h = \left(\sum\limits_{i=1}^n \vert y_i\vert^h\right)^{\frac{1}{h}}.$$

\noindent\begin{minipage}{\linewidth}
	\centering
	\begin{minipage}{0.48\linewidth}
		\begin{figure}[H]
	\includegraphics[width=0.7\linewidth]{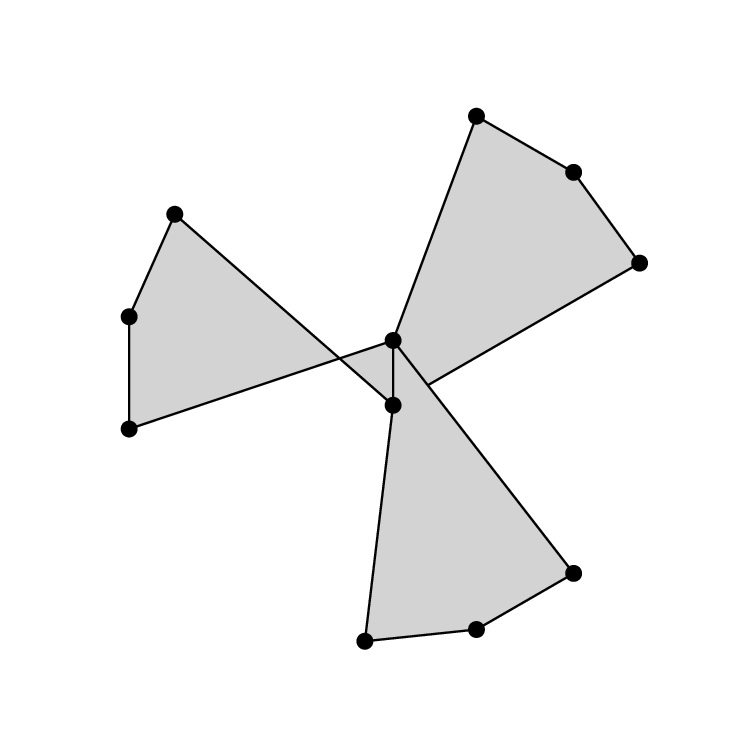}
			\caption{The sunflower hypergraph $\mathcal{S}(3,2,5)$.}
			\label{fig3}
		\end{figure}
	\end{minipage}
	\hspace{0.02\linewidth}
	\begin{minipage}{0.48\linewidth}
		\begin{figure}[H]
\includegraphics[width=0.7\linewidth]{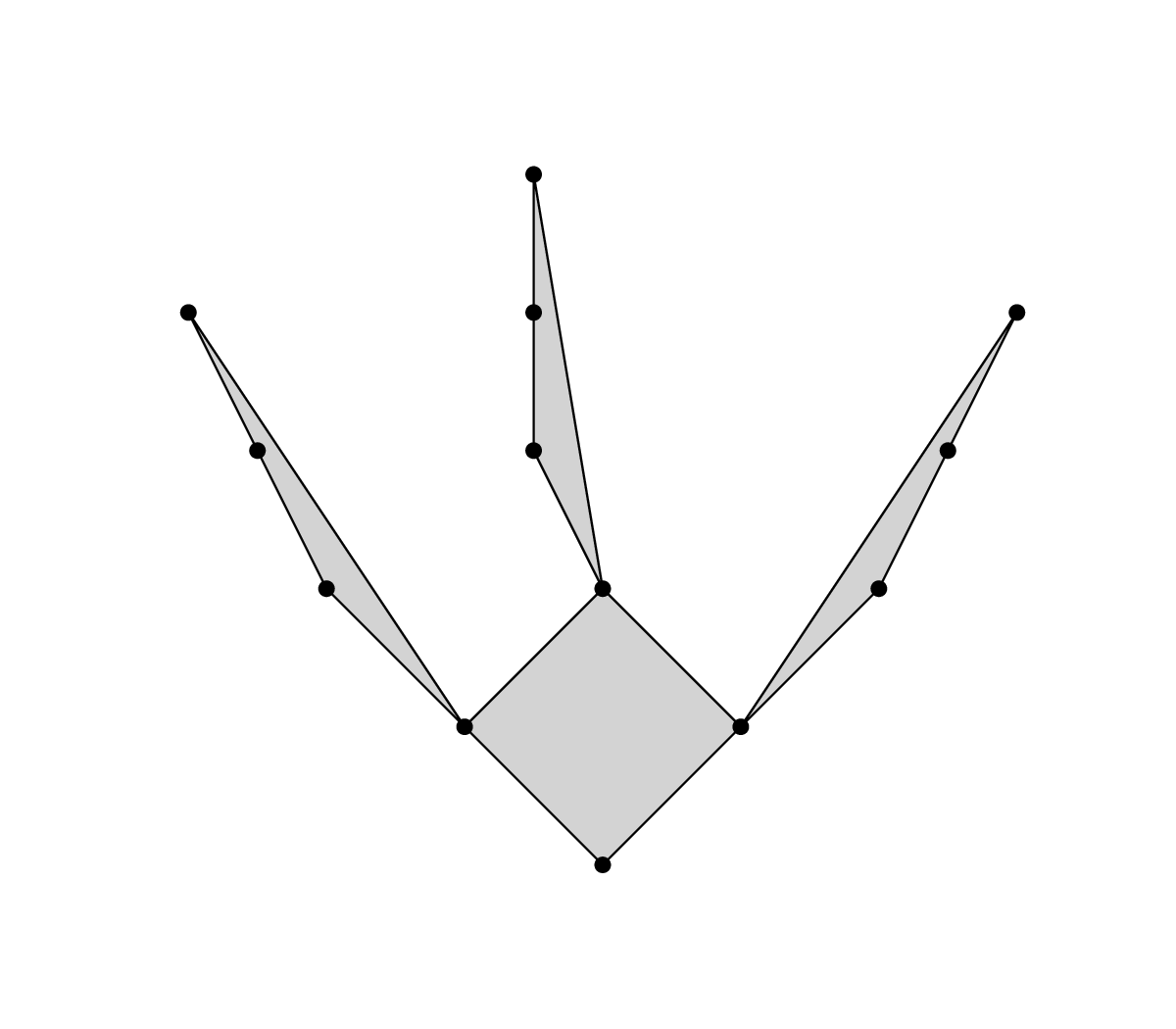}
			\caption{The $4$-uniform squid.}
			\label{fig4}
		\end{figure}
	\end{minipage}
\end{minipage}

\noindent Yue et al. \cite{yue2017adjacency} have obtained the largest eigenvalue of a sunflower (power hypergraph) and a squid (cored hypergraph).
The pictorial depiction of a $5$-uniform sunflower and a $4$-uniform squid are given in  Figure \ref{fig3} and Figure \ref{fig4}, respectively.

One of the most celebrated result in the field of spectral theory of graphs is about the largest eigenvalue of the matrix and its associated eigenvector, and this result is due to Perron and Frobenius.
The generalization of this result for non-negative hypermatrices has been obtained after a series of attempts. 

\begin{theorem}\cite{chang2008perron,yang2010further,yang2011some,friedland2013perron,qi2013symmetric}[Hypermatrix version of the theorem of Perron and Frobenius]
Let $\mathcal{M}$ be an order $h$ non-negative hypermatrix of dimension $n$. Then $\rho(\mathcal{M}) \geq 0$ is an eigenvalue of $\mathcal{M}$ corresponding to the non-zero eigenvector $\mathbf{y} \geq 0.$ 
    Furthermore, if $\mathcal{M}$ is weakly irreducible, then \\
    1. $\rho(\mathcal{M}) > 0$ is in $\sigma(\mathcal{M})$, and the associated eigenvector $\mathbf{y} > 0.$\\
    2. for any eigenvalue $\lambda_i$ of $\mathcal{M},$ it is true that $\vert \lambda_i \vert \leq \rho(\mathcal{M})$.\\
    3. an eigenvalue $\lambda_i$ corresponding to a non-negative (non-zero) eigenvector must be equal to $\rho(\mathcal{M}).$ 
Additionally, up to a multiplicative constant, the non-negative eigenvector (also known as the principal eigenvector or a Perron vector) is unique.
    \end{theorem}

For the hypermatrices associated with the hypergraphs, the existence and uniqueness conditions for the principal eigenvector associated with the largest eigenvalue is provided (with slightly different proof) by Pearson and Zhang \cite{pearson2014spectral}. 
It was proved \cite{fan2019eigenvectors} that $gm_{\mathcal{A}(\mathcal{H})}{(\rho)}$  =  $gm_{\mathcal{L}(\mathcal{H})}{(0)}$, where $\mathcal{L}(\mathcal{H})$ stands for the Laplacian hypermatrix \cite{qi2013h} of $\mathcal{H}$. 
The principal eigenvector is also considered as an eigenvector centrality, studying which is useful in ranking the vertices of the hypergraph. 
The principal ratio of $\mathcal{H}$ is defined to be the ratio 
$ \max_{ i \in [n]} y_i  / \min_{j \in [n]} y_j$, where $\mathbf{y}$ denotes the Perron eigenvector of $\mathcal{H}$. 
Numerous studies \cite{liu2016principal,li2018principal,cardoso2021principal,cooper2024principal,clark2023comparing,lia2024principal} on the principal eigenvector (mainly, the bounds for principal ratio) of the hypergraph were carried out, and consequently, the hypergraph's spectral radius.

\begin{theorem}\cite{deng2020geometry}
    Let $\mathcal{H}$ be a $h$-uniform regular hypergraph. Then the number of components in $\mathcal{H}$ is equal to the maximum number of linearly independent non-negative eigenvectors of $\mathcal{H}$ associated with $\rho(\mathcal{H})$. 
\end{theorem}

\noindent For a hypergraph $\mathcal{H}$, $p$-spectral radius is denoted by $\lambda^{(p)}(\mathcal{H}),$ and is defined as 
$$ \lambda_p(\mathcal{H}) = (h-1)! \max\left\{   \mathcal{A}(\mathcal{H}) \mathbf{y}^h : \left(\sum_{j=1}^n \vert y_j \vert \right)^{\frac{1}{p}} =1 \right\}.$$
It is a generalization of the Lagrangian of the hypergraph.
In addition to the study of the largest eigenvalue of the numerous hypermatrices associated with the hypergraph, much research has been done on the $p$-spectral radius \cite{chang2018computing,kang2015p,keevash2014spectral,nikiforov2014analytic,kang2014extremal,lu2018maximum,kang2019eigenvectors,nikiforov2019hoffman,zhou2020extremal,liu2022alpha}, which is not our focus here.

\begin{theorem}\cite{fan2022dimension}
   Suppose that $\mathcal{H} \in \prescript{}{m}{H}^{(h)}, m\geq 1$.
    \begin{itemize}
        \item If $\mathcal{H}$ is connected, then for all $\lambda \in \sigma(\mathcal{H})$ with $\vert \lambda \vert= \rho(\mathcal{H})$, the dimension of the eigen variety of $\mathcal{A}(\mathcal{H})$ associated with $\lambda$  is zero (i.e., it contains a finite number of eigenvectors). 
        \item If $\mathcal{H}$ is disconnected, then the dimension of the eigen variety associated with $\rho$ is greater than zero if and only if $\mathcal{H}$ consist of at least two components (of order greater than or equal to two) with spectral radius equal to $\rho$.
    \end{itemize}
\end{theorem}

Recently, Chen and Bu \cite{chen2024algebraic} have given the expression for the multiplicity of the largest eigenvalue of a hypertree of given size.

\begin{theorem}\cite{chen2024algebraic}
    If $\mathcal{T} \in \prescript{}{m}{H}^{(h)}$, then  $am_{\mathcal{T}}(\rho(\mathcal{T})) = h^{m(h-2)}.$
\end{theorem}

 Chen and et al. \cite{chen2024spectra} have found the multiplicity of the largest eigenvalue of a power of a simple graph. 
\begin{theorem}\cite{chen2024spectra}
  If $G \in \prescript{n'}{m}{H}^{(2)}$ is connected, then for $h \geq 3$,
   $$ am_{G^h}(\rho(G^h)) = h^{m(h-3)+n'-1}.$$
\end{theorem}

Bai and Lu \cite{bai2018bound} have characterized the hypergraph attaining the maximum value of the largest eigenvalue among all hypergraphs with $m$ hyperedges.

\begin{theorem}\cite{bai2018bound}
    If $\mathcal{H} \in \prescript{}{m}{H}^{(h)}$, then
    $$ \rho(\mathcal{H}) \leq g_h(m), \text{ where }  g_h(x) = p_{h-1}(p_h^{-1}(x)), \text{ and } p_h(x) = \binom{x}{h}.$$
    Upper bound equality is reached if $m = \binom{n_1}{h}$ for some $n_1 \geq h$ and $\mathcal{H} \cong K_{n_1}^{(h)} \cup \overline{K_{n-n_1}},$ where $n$ is the order of the hypergraph.
\end{theorem}
\noindent A similar type of the bound in the case of simple graphs was proved by Brualdi and Hoffman \cite{brualdi1985spectral}. 
They also conjectured a more general bound of the same type, which was later proven by Rowlinson \cite{rowlinson1988maximal}.

\begin{definition}
    Let $\mathcal{M}= (m_{i_1,\ldots,i_h})$ be an order $h$, dimension $n$ hypermatrix. Then, the $i^{th}$ row sum ($1 \leq i \leq n$) of $\mathcal{M}$, denoted by $S_i(\mathcal{M})$ is given by 
    $$ S_i(\mathcal{M}) = \sum\limits_{\substack{i_j \in [n]\\2 \leq j \leq h}} m_{ii_2\ldots i_h}.$$
\end{definition}

\begin{theorem}\cite{fan2015spectral,yang2010further}
    Let $\mathcal{M}= (m_{i_1,\ldots,i_h})$ be an order $h$, dimension $n$ hypermatrix with row sums $S_i$. Then, 
    $$ \min\limits_{i \in [n]}  S_i(\mathcal{M}) \leq \rho(\mathcal{M}) \leq \max\limits_{i \in [n]} S_i(\mathcal{M}).$$
    Furthermore, for a weakly irreducible hypermatrix $\mathcal{M},$  equality in the above is attained iff $S_i(\mathcal{M})=S_j(\mathcal{M})$ for every $i \neq j.$
\end{theorem}

From the above theorem it is direct that the largest eigenvalue of a hypergraph is  bounded  above by $\Delta$  and bounded below by $\delta$. 

\begin{theorem}\cite{you2019sharp}
    For $h\geq 3$, if $\mathcal{H}\in \prescript{n}{}{H}^{(h)}$, then 
    $$ \min\limits_{e \in \mathcal{E}(\mathcal{H})} \min\limits_{ \{ u, v \} \subset e} \sqrt{d_u d_v} \leq \rho(\mathcal{H}) \leq   \max\limits_{e \in \mathcal{E}(\mathcal{H})} \max\limits_{ \{ u, v \} \subset e} \sqrt{d_u d_v}. $$
    Furthermore, for a connected $\mathcal{H}$,  regularity of the hypergraph is the necessity and sufficiency for the equality.
\end{theorem}

\begin{definition}
    Suppose that $\mathcal{M}_1$ and $\mathcal{M}_2$ are order $h$ hypermatrices of dimension $n$. Then $\mathcal{M}_1$ and $\mathcal{M}_2$ are called diagonally similar, if there exists an $n \times n$ diagonal, non-singular matrix $P$ such that $ P^{-(h-1)} \mathcal{M}_1 P = \mathcal{M}_2$.
\end{definition}

\begin{theorem}\cite{shao2013general}
    Diagonally similar hypermatrices (order $h$, dimension $n$) have same spectra.
\end{theorem}

\begin{theorem}\cite{yang2011some}
    Let $\mathcal{M}_1$ and $\mathcal{M}_2$  be order $h$, dimension $n$ hypermatrices satisfying $\vert \mathcal{M}_1 \vert \leq \mathcal{M}_2$ and $\mathcal{M}_2$ is weakly irreducible. If $\lambda \in \sigma(\mathcal{M}_1$), then $\vert \lambda \vert \leq \rho(\mathcal{M}_2).$
\end{theorem}

\begin{definition}\label{blow_up}
Given a hypergraph $\mathcal{H}_0  \in \prescript{n-1}{}{H}^{(h-1)}$,
 the hypergraph $\mathcal{H} \in  \prescript{n}{}{H}^{(h)}$  obtained from $\mathcal{H}_0$ by including a vertex $v \notin \mathcal{V}(\mathcal{H}_0)$ to each hyperedge of $\mathcal{H}_0$ is called the blow-up of the hypergraph $\mathcal{H}_0$. 
\end{definition}

\begin{theorem}\cite{yuan2015some}
    If  $d_1 \geq  \cdots \geq d_n$ denotes the degree sequence of the hypergraph  $\mathcal{H} \in  \prescript{n}{}{H}^{(h)}$, then 
     $\rho(\mathcal{H})$ is upper bounded by   $$d_1^\frac{1}{h}d_2^{1-\frac{1}{h}},$$ and the necessary and sufficient condition for equality is either $\mathcal{H}$ is regular or isomorphic to a blow-up of a regular one. 
\end{theorem}

Later, the blow-up of the hypergraph is generalized \cite{xu2023extremal}, where each vertex of the hypergraph is replaced with a class (each of cardinality greater than one and they sum to $n$) of vertices. 
Among all possible blow-ups of the hypergraph (of fixed order), Xu, Hu and Wang \cite{xu2023extremal} have characterized the hypergraphs that attains the largest and least spectral radius.


\begin{theorem}\cite{nikiforov2014analytic}
    If $\mathcal{H} \in  \prescript{n}{m}{H}^{(h)}$, then $$ \rho(\mathcal{H}) \leq \frac{h}{\sqrt[h]{h!} } m^{\frac{h-1}{h}}.$$
    Moreover, if $\mathcal{H}$ is $h$-partite, then $\rho(\mathcal{H}) \leq m^{\frac{h-1}{h}},$ with necessary and sufficient condition is $\mathcal{H}$ being a complete $h$-partite.
\end{theorem}

Nikiforov \cite{nikiforov2014analytic} conjectured the following upper bound and it is proved by Kang, Liu and Shan \cite{kang2018sharp}.
\begin{theorem}\cite{kang2018sharp,kang2018spectral}
     If $\mathcal{H}\in  \prescript{n}{m}{H}^{(h)}$, then
    $$ \rho(\mathcal{H}) \geq \left( \frac{1}{n} \sum\limits_{u \in \mathcal{V}(\mathcal{H})} d_u^{\frac{h}{h-1}}  \right)^{\frac{h-1}{h}}. $$
 Furthermore, for a connected $\mathcal{H}$, regularity of the hypergraph is the necessity and sufficiency for the equality..
\end{theorem}

The following bounds are due to the bound for an irregularity measure $\rho-d_{avg}$.

\begin{theorem}\cite{liu2018irregularity}
    Let $\mathcal{H}\in  \prescript{n}{m}{H}^{(h)}$. If $d_{avg}=\frac{hm}{n}$ denotes the average degree of the hypergraph, then 
    $$ d_{avg} + \left(\frac{\sqrt[h]{h!}}{h^h} \right)^{\frac{1}{h-1}} \frac{h-1}{\sqrt[r]{m}} D_1 \leq \rho(\mathcal{H}) \leq d_{avg} + \left(\frac{D_2}{2} \right)^{\frac{h-1}{h}} \frac{h}{\sqrt[h]{h!}},$$
    where $D_1= \frac{1}{n} \sum\limits_{u \in \mathcal{V}(\mathcal{H})} d_u^{\frac{h}{h-1}} - d_{avg}^{\frac{h}{h-1}} $ and $ D_2 = \sum\limits_{u \in \mathcal{V}(\mathcal{H})} \left\vert d_u - d_{avg} \right\vert$.
\end{theorem}

\begin{theorem}\cite{lv2020sharp,chen2017spectral}
    Let $d_1 \geq \cdots \geq d_n$ be the degree sequence of $\mathcal{H} \in  \prescript{n}{m}{H}^{(h)}$. If  $P= \frac{1}{h-1}\binom{n-2}{h-2}$, then 
    $$ \rho(\mathcal{H}) \leq \min\limits_{1 \leq j \leq n} T_j,$$ where $T_1 = d_1$, and for $2 \leq j \leq n$, 
    $$ T_j = \frac{1}{2} \left( d_j - P + \sqrt{(d_j+P)^2 + 4P\sum\limits_{i=1}^{j-1} (d_i - d_j)} \right).$$
 Furthermore, for a connected $\mathcal{H}$,  regularity of the hypergraph is the necessity and sufficiency for the equality.
\end{theorem}

The generalization of the above result for the average $q$ degrees of the hypergraph has been carried out by the same authors \cite{lv2020general}.

Nikiforov \cite{nikiforov2014analytic} conjectured that the spectral radius of the shadow of a hypergraph $\mathcal{H}$ is at most $h-1$ times the spectral radius of $\mathcal{H}$. This conjecture is disproved by Liu, Kang and Bai \cite{liu2019bounds}, and they have obtained the following bounds for the spectral radius in terms of degree and co-degree of the vertices.
\begin{theorem}\cite{liu2019bounds}
If $d_{u_1u_2}$ is the co-degree of the vertices $u_1,u_2$ in $\mathcal{H} \in \prescript{n}{}{H}^{(h)}$, then 
$$ \rho(\mathcal{H}) \geq \frac{2 \sum_{\{u_1,u_2 \} \subseteq \mathcal{V}(\mathcal{H}), u_1 \neq u_2 } d_{u_1u_2} (d_{u_1} d_{u_2})^{\frac{1}{h-1}} }{ \left( \sum_{u \in \mathcal{V}(\mathcal{H})}d_u^{\frac{h}{h-1}} \right)^{\frac{2}{h}} (h-1)n^{1-\frac{2}{h}} }. $$
Furthermore for connected $\mathcal{H}$ with $h \geq 3$,  regularity of the hypergraph is the necessity and sufficiency for the equality.
\end{theorem}

The co-degree matrix $\mathcal{C}(\mathcal{H})$ of $\mathcal{H}$  is a square matrix with zero diagonal, and for $u_1 \neq u_2$, $(u_1,u_2)$-th entry is $d_{u_1u_2}$. 
In this result, the authors have obtained the relation between $\rho(\mathcal{H})$ and $\rho(\mathcal{C}(\mathcal{H}))$.
\begin{theorem}\cite{liu2019bounds}
    For $\mathcal{H} \in \prescript{n}{}{H}^{(h)}$, let $\mathcal{C}(\mathcal{H})$ denotes the co-degree matrix $\mathcal{H}$.  Then
    $$ \rho(\mathcal{H}) \geq \frac{1}{h-1} \left(\frac{2}{n} \right)^{\frac{h-2}{h}} \rho(\mathcal{C}(\mathcal{H})). $$
\end{theorem}
Coloring of a hypergraph $\mathcal{H}$ is ascribing the vertices of $\mathcal{H}$ with colors in such a way that no hyperedge is monochromatic. 
If a hypergraph $\mathcal{H}$ can be colored with $t$ colors, then $\mathcal{H}$ is $t$-colorable. 
$\chi(\mathcal{H})$ denotes the chromatic number of $\mathcal{H}$, which is the least number of colors needed to color all the vertices.

\begin{theorem}\cite{cooper2012spectra}
    Let $\Delta$ be the maximum degree, $d_{avg}$ be the average degree, and $\chi(\mathcal{H})$ be the chromatic number of $\mathcal{H} \in \prescript{}{}{H}^{(h)}$.  Then \\
    1. $\rho(\mathcal{H}) \geq \chi(\mathcal{H})-1$;\\
    2. $d_{avg} \leq \rho(\mathcal{H}) \leq \Delta, $ and equality is attained for a regular hypergraph.
\end{theorem}

\begin{definition}\cite{lin2016sharp}
    Let $\mathcal{M}= (m_{i_1,\ldots,i_h})$ be an order $h$, dimension $n$ hypermatrix and $s_i(\mathcal{M}) (> 0)$ be the row sum corresponding to the index $i$. Then the average $2$-row sum of the hypermatrix $\mathcal{M}$ is denoted by $S_i^{(2)}(\mathcal{M})$ and is defined as 
    \begin{equation*}
        S_i^{(2)}(\mathcal{M}) = \frac{\sum_{\substack{i_j \in [n]\\ 2 \leq j \leq h}} m_{ii_2\ldots i_h} S_{i_2}(\mathcal{M})\cdots S_{i_h}(\mathcal{M})}{(S_i(\mathcal{M}))^{h-1}}.
    \end{equation*}
\end{definition}

In the following result, using the average $2$-degree of the hypergraph the bound has been obtained for the largest eigenvalue.

\begin{theorem}\cite{lin2016sharp,bu2017brauer}
    Let $\mathcal{M}$ be an order $h$ hypermatrix of dimension $n$. If $S_i^{(2)}(\mathcal{M})$ denotes the average $i$-{th}  $2$-row sum of the hypermatrix $\mathcal{M},$ then 
    \begin{equation*}
        \min\limits_{1 \leq i \leq n} S_i^{(2)}(\mathcal{M})  \leq \rho(\mathcal{M}) \leq \max\limits_{1 \leq i \leq n} S_i^{(2)}(\mathcal{M})
    \end{equation*}
    with equality holds if and only if $S_1^{(2)}(\mathcal{M})=\cdots=S_n^{(2)}(\mathcal{M}).$ 
\end{theorem}

\begin{definition}
Suppose $d_1,\ldots,d_n$ denotes the vertex degrees of $\mathcal{H}$. The average $2$-degree of a vertex $i$ is $t_i= \frac{\sum_{\substack{e \in \mathcal{E}_i}}  d_{i_2}\cdots d_{i_h}}{d_i^{h-1}},$  where \\ $e=\{ i, i_2,\ldots,i_h\}.$
\end{definition}

\begin{theorem}\cite{lin2016sharp,bu2017brauer}
    Let $t_1 \geq  \cdots \geq t_n$ be the  average $2$-degrees of $\mathcal{H}\in \prescript{n}{}{H}^{(h)}$. Then 
    $$\rho(\mathcal{H}) \leq t_1^{\frac{1}{h}}t_2^{1-\frac{1}{h}}, $$
    with necessity and sufficiency for the equality is when every vertex has same average $2$-degrees.
\end{theorem}

The generalization of the definition of the average $2$ degrees has been done to average $q$ degrees, for some $q \geq 2$. 
Furthermore, the similar type of upper bound holds \cite{fang2023some} for the largest eigenvalue of the hypergraph using the average $q$ degrees.
Fang, Huang, and You \cite{fang2023some} have also discussed the sharpness of this bound.
The Nordhaus-Gaddum type of bounds for the largest eigenvalue of a hypergraph are obtained in \cite{si2017spectral}.

Lu and Man \cite{lu2016connected} have used a novel technique to determine (or to bound) the largest eigenvalue of a $h$-uniform hypergraphs using some labeling techniques. 
They have characterized all hypergraphs (uniform) with with largest eigenvalue at most (and, equal to) $\sqrt[r]{4}$.
The ordering of these hypergraphs with respect to their spectral radius was carried out by Li and Chang \cite{li2020effect}.

\begin{definition}\cite{lu2016connected}
    Given a hypergraph $\mathcal{H}$ of order $n$ and size $m$, a weighted incidence matrix $B$ of $\mathcal{H}$ is an $n \times m$ matrix, whose rows correspond to the vertices, columns correspond to the hyperedges of the hypergraph, and (for any vertex $u$ and hyperedge $e$  of $\mathcal{H}$) $B(u,e) >0$ if $u \in e$ and $B(u,e)=0$ otherwise.
\end{definition}

\begin{definition}\cite{lu2016connected}
     Given a $h$-uniform hypergraph $\mathcal{H}$, if there exist a incidence (weighted) matrix $B$ of $\mathcal{H}$ satisfying the following two conditions, then we say that $\mathcal{H}$ is $\alpha$-normal. That is
     \begin{itemize}
        \item[1.] $\sum\limits_{e \in \mathcal{E}_u} B(u,e)=1,$ (where $\mathcal{E}_u = \{ e \in \mathcal{E} \vert  u \in e\}$);
        \item[2.] $\prod\limits_{u \in e} B(u,e)=\alpha.$
    Furthermore,  if for any cycle (if exists) $u_0 e_1 u_1 e_2 \ldots e_t \\ u_0, $ $\prod\limits_{i=1}^t \frac{B(u_i-1,e_i)}{B(u_i,e_i)}=1$, then we say that $\mathcal{H}$ is consistent. 
    \end{itemize}
\end{definition}

\begin{theorem}\cite{lu2016connected}
    For a connected $\mathcal{H} \in H^{(h)}$, we have  $\rho(\mathcal{H})=\mu$ iff $\mathcal{H}$ is consistent and $\mu^{-r}$-normal.
\end{theorem}

\begin{definition}\cite{lu2016connected}
    Given a $h$-uniform hypergraph $\mathcal{H}$, if there exist a incidence (weighted) matrix $B$ of $\mathcal{H}$ satisfying the following two conditions, then we say that $\mathcal{H}$ is $\alpha$-subnormal (resp. strictly $\alpha$-supernormal). That are
    \begin{itemize}
        \item[1.] $\sum\limits_{e \in \mathcal{E}_u} B(u,e) \leq 1,$ (resp. $\sum\limits_{e \in \mathcal{E}_u} B(u,e) \geq 1,$);
        \item[2.] $\prod\limits_{u \in e} B(u,e) \geq \alpha$ (resp. $\prod\limits_{u \in e} B(u,e) \leq \alpha$).
    \end{itemize}
    Moreover, if $\mathcal{H}$ is $\alpha$-subnormal (resp.  $\alpha$-supernormal) but not $\alpha$-normal, then we say that $\mathcal{H}$ is strict $\alpha$-subnormal (resp. strict $\alpha$-supernormal).
\end{definition}

\begin{theorem}\cite{lu2016connected}
    For a connected  $\alpha$-subnormal hypergraph $\mathcal{H} \in H^{(h)}$, we have $$\rho(\mathcal{H}) \leq \alpha^{-\frac{1}{h}}.$$
 For a connected  $\alpha$-supernormal hypergraph $\mathcal{H} \in H^{(h)}$, we have  $$\rho(\mathcal{H}) \geq \alpha^{-\frac{1}{h}}.$$
\end{theorem}

Li, Zhou and Bu \cite{li2018principal} have got an upper bound for the largest eigenvalue of a connected non-regular hypergraphs of given diameter, size and order. 
The study of $\Delta-\rho$, another kind of irregularity measure, makes up for the following bounds.
 A similar type of bound is obtained for simple irregular graphs by Cioab{\u{a}}, Gregory and Nikiforov \cite{cioabua2007extreme}.

\begin{theorem}\cite{li2018principal}
    Let $\mathcal{H} \in \prescript{n}{m}{H}^{(h)}$, $h \geq 3$  be connected and irregular (not regular) with diameter $D$. If $\Delta$ denotes the maximum degree, then
    $$ \rho(\mathcal{H}) < \Delta + \frac{\Delta( \Delta n-mr)}{2mD(h-1)(\Delta n-mr)+\Delta n}. $$
\end{theorem}

Some other works on the irregularity measure of hypergraphs were carried out in \cite{zhou2015spectral}. 
The ordering of hypertrees and unicyclic hypergraphs with respect to the irregularity measure $\Delta-\rho$ was carried out by Lin, Yu and Zhou \cite{lin2023irregularity}.

\begin{definition}\cite{kang2018spectral}
    Given a hypergraph $\mathcal{H}=({V},{E})$ with uniformity $k$, and for any $h \geq k,$ $ 1 \leq s \leq \lfloor \frac{h}{k}  \rfloor.$ Define pairwise disjoint set $V_u$ for each $u \in \mathcal{V}$ and $V_e$  for each $e \in \mathcal{E}.$ Now, the generalized power hypergraph $\mathcal{H}=H^{h,s}$ is a hypergraph with uniformity $h$ on vertex set $\mathcal{V}_u= \left( \bigcup\limits_{ u \in V} \mathcal{V}_u \right) \cup \left( \bigcup\limits_{e \in {E}} \mathcal{V}_e \right)$ and with  $\mathcal{E}(\mathcal{H}) = \left\{ \left( \bigcup\limits_{u \in e} \mathcal{V}_u \right) \cup \mathcal{V}_e : e \in E  \right\}.$ 
\end{definition}

In \cite{fan2016h} Khan, Fan and Tan have shown that, the spectrum of a simple graph $G$ (as a set) is contained in the spectrum of $G^{h,s}$.

\begin{theorem}\cite{kang2018spectral}
    Let $\mathcal{G} \in H^{(k)}$ and $\mathcal{H}=\mathcal{G}^{h,s}$, $h \geq k \geq 2,$ $ 1 \leq s \leq \lfloor \frac{h}{k}  \rfloor$ be a generalized power hypergraph. If $ 0 < \lambda \in  \sigma(\mathcal{G})$ is corresponding to a non-negative eigenvector, then $\lambda^{\frac{ks}{h}} \in \sigma(\mathcal{H}).$ Moreover, for a connected $\mathcal{G},$ $\rho(H) = \rho(\mathcal{H})^{\frac{ks}{h}}$.
\end{theorem}

\begin{definition}\cite{li2016extremal}
 Given a hypergraph $\mathcal{H}=(\mathcal{V},\mathcal{E})$ and a vertex $u \in \mathcal{V}$, let $e_1,\ldots,e_t \in \mathcal{E}$ be the hyperedges of $\mathcal{H}$ such that $u \notin e_i,$ for any $i, 1\leq i \leq t.$ Let $u_i \in e_i$ and $e_i'=(e_i \setminus u_i)\cup \{ u \}$ for each $i.$ Also, let $\mathcal{H}'=(\mathcal{V},\mathcal{E}')$ be a new hypergraph with the edge set $\mathcal{E}' = (\mathcal{E}\setminus \{ e_1,\ldots,e_t \}) \cup \{ e_1',\ldots,e_t' \}.$ Then we say that the hypergraph $\mathcal{H}'$ is obtained from moving the hyperedges $(e_1,\ldots,e_h)$ from $(v_1,\ldots,v_h)$ to $u$.
\end{definition}

\begin{lemma}\cite{li2016extremal}
    Given a connected hypergraph $\mathcal{H} \in H^{(r)}$, let $\mathcal{H}'$ be obtained from $\mathcal{H}$ by moving the hyperedges $(e_1, \ldots,e_t)$ from $(v_1,\ldots, v_t)$ to $u$  (where $t \geq 1$), where $\mathcal{H}'$ is also a simple hypergraph. If $y$ be the Perron (eigen) vector of $\mathcal{A}(\mathcal{H})$ corresponding to $\rho(\mathcal{H}),$ and $y_u \geq \max\limits_{1 \leq i \leq t} y_{u_i}$  then $\rho(\mathcal{H}) < \rho(\mathcal{H}').$   
    \end{lemma}

\begin{lemma}\cite{fan2016maximizing}
    If $\mathcal{H}$ is a maximizing (spectral radius) $h$-uniform hypergraph among all connected $h$-uniform hypergraph of fixed size, then $\mathcal{H}$ contains a vertex of full degree. 
\end{lemma}

Any connected $h$-uniform hypergraph with $m$ hyperedges on $m(h-1)-1+l$ vertices is called a $l$-cyclic hypergraph.
Authors of the article \cite{fan2016maximizing}, have also given the bound for the largest eigenvalue of the hypergraphs among the family of (linear or power) unicyclic and (linear or power) bicyclic hypergraphs.
Ouyang, Qi and Yuan \cite{ouyang2017first} have characterized the first three bicyclic, five unicyclic hypergraphs attaining maximum value of the largest eigenvalue among the family of unicyclic and bicyclic uniform hypergraphs, respectively.
It was found that the unicyclic uniform hypergraph attaining the maximum value of the largest eigenvalue is not linear.
Ding, Fan and Wan \cite{ding2020linear} have characterized the hypergraph with the maximum, the second maximum and the third maximum value of the largest eigenvalue in the class of all unicyclic linear hypergraphs.
Shan, Wang and Wang \cite{shan2021smallest} have characterized the hypergraph attaining the minimum value of the largest eigenvalue in the class of all bicyclic hypergraphs of given size.
Zheng et al. \cite{zheng2024spectral} have characterized first eight hypergraphs of size $m > 20,$  attaining the maximum value of the largest spectral radius among the class of all tricyclic uniform hypergraphs. 

Xiao and Wang \cite{xiao2020effect} studied the effect of subdividing a hyperedge of the hyperpath on the largest eigenvalue of the hypergraph. 
They proved that the largest eigenvalue of the hypergraph increases (resp. decreases) if the subdividing hyperedge is a part of the pendent path (resp. the internal path \cite{xiao2020effect}) of the hypergraph.

\begin{theorem}\cite{lin2016upper}
    If $\mathcal{H} \in \prescript{n}{m}{H}^{(h)}$, then 
    \begin{equation*}
        \rho(\mathcal{H}) \leq \frac{h m^{1-\frac{1}{h}} \binom{n}{h}^{\frac{1}{h}}}{n}.
    \end{equation*}
    Equality is reached if and only if $\mathcal{H} \cong \overline{\mathcal{K}_n}$  or $\mathcal{H} \cong {\mathcal{K}_n^{(h)}}$.
\end{theorem}

\begin{theorem}\cite{lin2016upper}
    If  $\mathcal{H} \in \prescript{n}{m}{H}^{(h)}$ is linear with minimum and maximum vertex degrees  $\delta \geq 1$ and $\Delta$, respectively, then
    \begin{equation*}
    \rho(\mathcal{H}) \leq \left(\frac{hm+\Delta(\delta(h-1)-1)-\delta(n-1))}{h-1}\right)^{\frac{h-1}{h}}.
    \end{equation*}
  The attainment of the equality happens if and only if $\mathcal{H}$ is either $2$-uniform and regular, or it is $2$-uniform and $ \mathcal{H} \cong t K_2 \cup S_{n-2t}, 0 \leq t \leq \left\lfloor \frac{n-3}{2} \right\rfloor,$ or $3.$ $ K_{\Delta+1} \cup \mathcal{H}',$ where $\mathcal{H}'$ is regular (simple graph) with regularity $\delta < \Delta.$
\end{theorem}

\begin{theorem}\cite{wan2022distribution}
    If $\Delta$ is the maximum vertex degree of $\mathcal{H} \in H^{(h)}$, then 
    $$ \Delta^{\frac{1}{h}} \leq \rho(\mathcal{H}) < \frac{h}{h-1} ((\Delta-1)(h-1))^{\frac{1}{h}}.$$
\end{theorem}

Wang and Yuan \cite{wang2020uniform} have characterized the hypertree that attains the minimum value of the largest eigenvalue among the class of  hypertrees of fixed maximum degree. 
Wang \cite{wang2022minimum} has characterized the hypertree attaining the minimum value of the largest eigenvalue in the class of all hypergraphs with size fixed and having two vertices of degree $\Delta$. He has also characterized (using a simpler method) the hypertree attaining the minimum value of the largest eigenvalue in the class of all hypertrees of fixed size.

A vertex subset $D \subseteq \mathcal{V}$ of $\mathcal{H}=(\mathcal{V},\mathcal{E})$ is called a dominating set, if for every   $ u \in \mathcal{V} \setminus D$, the existence of a hyperedge $e \in \mathcal{E}$ such that $e \cap D \neq \emptyset$ is guaranteed. 
We call a dominating set $D$ as a minimal dominating set if  for any $u \in D$, the set  $D \setminus \{u\}$ does not dominate $\mathcal{H}$.
A dominating set that is minimal and of minimum cardinality is a minimum dominating set, and the cardinality is defined to be the domination number, which we denote by $\gamma(\mathcal{H})$ (or simply $\gamma$).
In the following theorem, a bound for the largest eigenvalue of hypergraphs using the domination number obtained and the equality case was also discussed. 

\begin{theorem}\cite{kang2017spectral,kang2021spectral}
    Let $\mathcal{H} \in \prescript{n}{}{H}^{(h)}$ be linear with domination number $\gamma$. Then 
    $$ \rho(\mathcal{H}) \geq \left( \left\lceil \frac{ \lceil n / \gamma -1 \rceil}{h-1} \right\rceil \right)^{\frac{1}{h}}.$$
    Equality is attained if and only if $\mathcal{H} = \mathcal{H}' \cup \mathcal{H}''$, where$ \mathcal{H}'\cong \mathcal{S}_{m'}^{(h)}$ with $m' = \left\lceil \frac{ \lceil n / \gamma -1 \rceil}{h-1} \right\rceil$, and $\mathcal{H}''$ is a $h$-uniform hypergraph having domination number $\gamma-1$ whose largest eigenvalue is upper bounded by $\left( \left\lceil \frac{ \lceil n / \gamma -1 \rceil}{h-1} \right\rceil \right)^{\frac{1}{h}}$.
\end{theorem}

Kannan, Shaked and Berman \cite{rajesh2016weakly} have shown that the spectral radius of a hypermatrix is always greater than the spectral radius of its (proper) principal sub-hypermatrix when it is weakly irreducible.

\begin{theorem}\cite{rajesh2016weakly}
   For an order $h$ weakly irreducible hypermatrix $\mathcal{M}$ of dimension $n$,if  $\mathcal{M}[I], I \subset [n]$ denotes the principal sub-hypermatrix of $\mathcal{M}$, then $\rho(\mathcal{M}[I]) < {\rho(\mathcal{M})}.$ 
\end{theorem}

Given a hypergraph $\mathcal{H}$ and a vertex $u$ it, if $\mathcal{H}-u$ denotes the sub-hypergraph (induced) of $\mathcal{H}$ induced by $\mathcal{V} \setminus \{ u \}.$
In the following the theorems, the bounds for the largest eigenvalue of the vertex deleted sub-hypergraph of the hypergraph is obtained.

\begin{theorem}\cite{lin2023largest}
    For a $\mathcal{H} \in H^{(h)}$, let  $\mathbf{y}$ be the non-negative eigenvector corresponding to $\rho(\mathcal{H})$ with $|| \mathbf{y} ||_{h} =1.$ Let $u \in \mathcal{V}(\mathcal{H})$ be fixed and  $\mathbf{z}$ be the restriction of $\mathbf{y}$ on $\mathcal{V}(\mathcal{H}) \setminus \{u \}.$
    Then  
    \begin{equation*}
        (1-ry_u^h) \rho(\mathcal{H}) \leq \rho(\mathcal{H} - u) \leq \rho(\mathcal{H}).
    \end{equation*}
Moreover, if $y_u^h \neq 1,$ then $\rho \left( \mathcal{H} - u \right) \geq \left( \frac{1-ry_u^h}{1-y_u^h}\right),$ and for a connected hypergraph $\mathcal{H}$, attainment of equality happens in the lower bound iff  $\mathbf{z}$ is the eigenvector of $\mathcal{H} - u$ corresponding to $\rho(\mathcal{H} - u).$
\end{theorem}

\begin{theorem}\cite{nikiforov2014analytic, lin2023largest}
    Let $\mathcal{H} \in \prescript{n}{}{H}^{(r)}$ with vertex set $[n]$.
    For any vertex $u \in [n],$ 
    \begin{equation*}
        \rho(\mathcal{H} - u) \geq \rho(\mathcal{H}) - \sqrt[h-1]{\frac{d_{\mathcal{H}}(u)}{\rho(\mathcal{H})}},
    \end{equation*}
    with equality if and only if $u$ is not an isolated vertex, and $\mathcal{H} - u$ is regular.
\end{theorem}

In addition to analyzing the largest eigenvalue of the hypergraph, research is also being done on other eigenvalues of the hypergraph, albeit to a lesser extent.
If a hypergraph is regular, then the largest eigenvalue will be its highest degree. So here comes the motivation for studying the second largest eigenvalue of the hypergraphs. In this context, Alon has proved \cite{nilli1991second} a famous lower bound (\emph{Alon-Boppana bound}) for the second largest eigenvalue of the regular graphs. The Alon-Boppana type bound for the eigenvalue (second largest) in the case of a regular uniform hypergraph was proposed by Li and Mohar \cite{li2019first}. 
 Later,  R{\"a}ty and Tomon  \cite{raty2024bisection} improved the above bound of Li and Mohar and also studied minimum bisection problem for hypergraphs. 
 The generalization of the above theorem for hypergraphs using the adjacency matrix also attracted much attention \cite{cioabua2022spectrum,lato2023spectral,song2023hypergraph}.
Fan, Zhu and Wang \cite{fan2020least} have characterized the hypergraph that attains least $H$-eigenvalue among the hypergraph with cut vertices.
 Kenter \cite{kenter2014necessary} has given the bounds for the smallest (the second smallest) eigenvalue with the average degree $d_{avg}$ of the hypergraph for a $2$-colorable $h$-uniform ($h=3,4$ and $5$) hypergraph.

\subsection{Spectral Tur\'an type problems}

A survey on the results of spectral analogues of the theory of extremal graphs can be found in \cite{nikiforov2011some}.
In literature, the spectral version of the hypergraph Tur\'an type problems were also discussed, which was initiated in \cite{hou2021spectral}.
 If a hypergraph does not contain $\Gamma$ as an induced sub-hypergraph, then we call it a $\Gamma$-free hypergraph.
If $\mathcal{G}$ denotes a family of hypergraphs and suppose $\mathcal{H}$ does not contain $G$ as an induced sub-hypergraph, for any $G \in \mathcal{F}$, then we say that  $\mathcal{H}$ is  $\mathcal{G}$-free. 

Let $spex_h(n, \mathcal{G})$ (resp. $spex_h^{lin}(n, \mathcal{G})$) denote the maximum value of spectral radius of $h$-uniform (resp.  $h$-uniform linear) $\mathcal{G}$-free  hypergraph on $n$ vertices.
Given a simple graph $G$ and a  hypergraph $\mathcal{H}$, we say that $\mathcal{H}$ is a Berge-$G$ hypergraph if $\mathcal{V}(G) \subseteq \mathcal{V}(\mathcal{H}),$ and there is a bijection $\eta: \mathcal{E}(G) \longrightarrow \mathcal{E}(\mathcal{H})$ such that $e \subseteq \eta(e)$. 
Denote by $\mathcal{B}_h(G)$, the family of all $h$-uniform hypergraphs thatr are Berge-$G$.
She, Fan and Kang \cite{she2025spectral} have deduced that, for any $\delta > 0$ and large enough $n$, 
$$n^{1-\delta} < spex_h^{lin}(n, \mathcal{B}_h(C_4) ) = \mathcal{O}(n). $$
Hou, Chang and Cooper \cite{hou2021spectral} have shown that, if $C_4$ denotes a cycle on $4$ vertices, then 
$$ spex_h^{lin}(n, \mathcal{B}_h(C_4) ) \leq \frac{\sqrt{n}}{h-1} + \mathcal{O}(1).$$
A $h$-uniform $Fan$, $F^h$  is a $h$-uniform hypergraph  obtained from a $h$-uniform hyperstar of size $h$, by adding a new hyperedge which is formed by including exactly one  pendent vertex from each of the $h$ hyperedges. They \cite{hou2021spectral} have also shown that, for $n \equiv 0 \mod{h},$
$$ spex_h^{lin}(n, F^h ) = \frac{n}{h}.$$
Let $K_{h+1}^h$ be the $h$ power hypergraph of the complete graph $K_{h+1}$.
 Gao, Chang and Hou \cite{gao2022spectral} have shown that, 
$$ spex_h^{lin}(n, K_{h+1}^h) \leq \frac{n}{h},$$ 
and the  attainment of the equality happens if and only if $h | n$ and the hypergraph is a balanced $h$-uniform, $h$-partite such that any two vertices from the distinct partite sets is contained in one hyperedge only.

Suppose that $G$ is a simple graph with chromatic number $t, ~ t > h\geq 2$. If $G^h$ denote an $h$-power hypergraph of the graph $G$, then \cite{she2022linear}
$$ spex_h^{lin}(n, G^h) \leq n \left( \frac{1}{h-1} \left( 1-\frac{1}{t-1}\right) + \mathcal{O}(1)\right).$$ 
A  bipartite complete graph with partition of the vertex set $\mathcal{V}=\mathcal{U}_1 \cup \mathcal{U}_2$, where $\vert \mathcal{U}_1 \vert =p $ and $\vert \mathcal{U}_2 \vert = q$  is denoted by $K_{p,q}$.
Zhu et al. \cite{zhu2024spectral} have shown that, 

$$spex_h^{lin}(n, \mathcal{B}_h(K_{2,q})) \leq \left(\left[ \frac{3(q-1)}{2} -\frac{q-3}{2(h-1)} \right] (n-1)\right)^{\frac{1}{2}}.$$
She, Fan and Kang \cite{she2025spectral} have proved the following:
\begin{itemize}
    \item For $ q \geq 2,$ if $4n(q-2) \geq (2h^2-4r+1)^2, $ then
    $$spex_h^{lin}(n, \mathcal{B}_h(K_{2,q})) \leq n^{\frac{1}{2}} \frac{(q-1)^{\frac{1}{2}}}{(h-1)} + n^{\frac{1}{4}} \frac{q^{\frac{1}{2}} (q-1)^{\frac{1}{4}}(2h^2-4r+1)^{\frac{1}{2}} }{h-1}.  $$
\item For $2 \leq p \leq q$,  if $\mathcal{F}= \mathcal{B}_h(K_{p,q}) \cup \mathcal{B}_h(C_3),$ then 
\begin{equation*}
    spex_h^{lin}(n, \mathcal{F}) \leq 
    \begin{cases}
        \frac{h-q+\sqrt{4(q-1)(n-1) + (h-q)^2}}{2(h-1)} & \text{ if } p=2, \\
        (p-2)+ \frac{p-1}{h-1} n^{1-\frac{2}{p}} + \frac{(q-p+1)^\frac{1}{p}}{h-1} n^{1-\frac{1}{p}} & \text{ if } p \geq 3.
    \end{cases}
\end{equation*}
\end{itemize}
For $n \geq h \geq 3$ and $q \geq p \geq 3,$  Zhou, Yuan and Chen \cite{zhou2025note} has proved the following upper bound. 
$$ spex_h^{lin}(n, \mathcal{B}_h(K_{p,q}) ) \leq \frac{n}{{h-1}} \sqrt{\frac{h(h-3)+4}{2(h-1)}} + \mathcal{O}(n^{\frac{1}{2}}). $$
Some more advancements in the spectral versions of Tur\'an type problems (mostly on $p$-spectral radius) can be found in \cite{nikiforov2017combinatorial,ni2024spectral,liu2024hypergraph,zheng2024spectralturan,cooper2024principal}.

If there exist a planar embedding for a graph simple graph $G$, whose all vertices lie in the boundary of the outer face, then $G$ is called an outer-planar graph.
A $3$-uniform hypergraph $\mathcal{H}$ is said to be an outer-planar \cite{zykov1974hypergraphs} if the shadow of $\mathcal{H}$ has the outer-planar embedding such that every hyperedge of $\mathcal{H}$ is the vertex set of an interior triangular face of the shadow.
Ellingham, Lu and Wang \cite{ellingham2022maximum} have characterized the $3$-uniform outer-planar hypergraph attaining the maximum value of the largest eigenvalue and they posed a conjecture (a bound) for $3$-uniform planar hypergraphs. 

\begin{theorem}\cite{ellingham2022maximum}
    Let $\mathcal{H} \in {n}{}{H}^{(3)}$ be an outer-planar hypergraph. For large enough $n$,
    $$ \rho(\mathcal{H}) \leq \rho(F_n^{(3)}),$$ where $F_n^{(3)}$ is a $3$-uniform hypergraph whose shadow is isomorphic to the join of a path of length $n-1$ with an isolated vertex.     
\end{theorem}

\subsection{Spectral radius of hypertrees}

In this section, we mainly focus on the results related to the eigenvalues of hypertrees and the characterizations of the extremal hypertrees.

It has been shown that \cite{clark2018adjacency} every  $h$-uniform hypertree is $h$-partite. 
Hence $h$-uniform hypertee is $h$-symmetric.

\begin{theorem}\cite{clark2018adjacency}
    Let $\sigma(\mathcal{T})$ be the spectrum of $\mathcal{T} \in T^{(h)}$ with $h\geq 3$ and $\zeta_h$ be the principal $h-th$ root of unity. Then $\sigma(\mathcal{T}) \subseteq \mathbb{R}[\zeta_h]$ if and only if $\mathcal{T}$ is a power hypertree.
\end{theorem}

In 1997, Hofmeister \cite{hofmeister1997two} has characterized trees that attains the largest and second largest spectral radius among the hypertrees of order $n$. Chang and Huang \cite{chang2003ordering} in 2003 have ordered first seven trees that attains the maximum spectral radius among all trees of order $n$. In 2006, Lin and Guo \cite{lin2006ordering}, ordered first thirteen trees that attains the maximum spectral radius among all trees of order $n$.

\begin{theorem}\cite{li2016extremal}
   Let $S_m^{(h)}$ and $P_m^{(h)}$, respectively be the $h$-uniform hyperstar and hyperpath with $m$ hyperedges. Then for any $\mathcal{T} \in T^{(h)}$,
    \begin{equation*}
        \rho(P_m^{(h)}) \leq \rho(\mathcal{T}) \leq \rho(S_m^{(h)}).
    \end{equation*}
\end{theorem}

The pictorial depiction of a $3$-uniform hyperpath and a  $3$-uniform hyperstar are presented in Figure \ref{fig1} and Figure \ref{fig2}, respectively.

\noindent\begin{minipage}{\linewidth}
	\centering
	\begin{minipage}{0.48\linewidth}
		\begin{figure}[H]
			\includegraphics[width=0.9\linewidth]{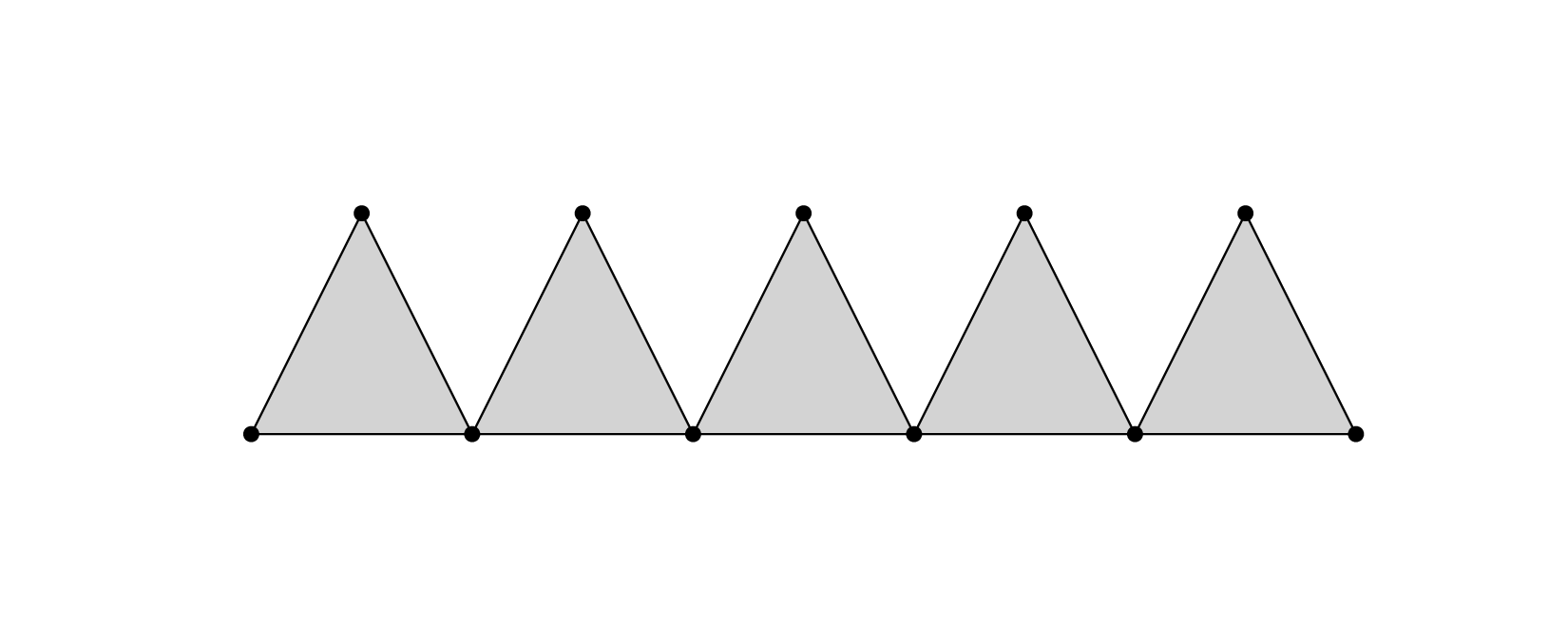}
			\caption{$3$-uniform hyperpath with $5$ hyperedges}
			\label{fig1}
		\end{figure}
	\end{minipage}
	\hspace{0.02\linewidth}
	\begin{minipage}{0.48\linewidth}
		\begin{figure}[H]
			\includegraphics[width=0.9\linewidth]{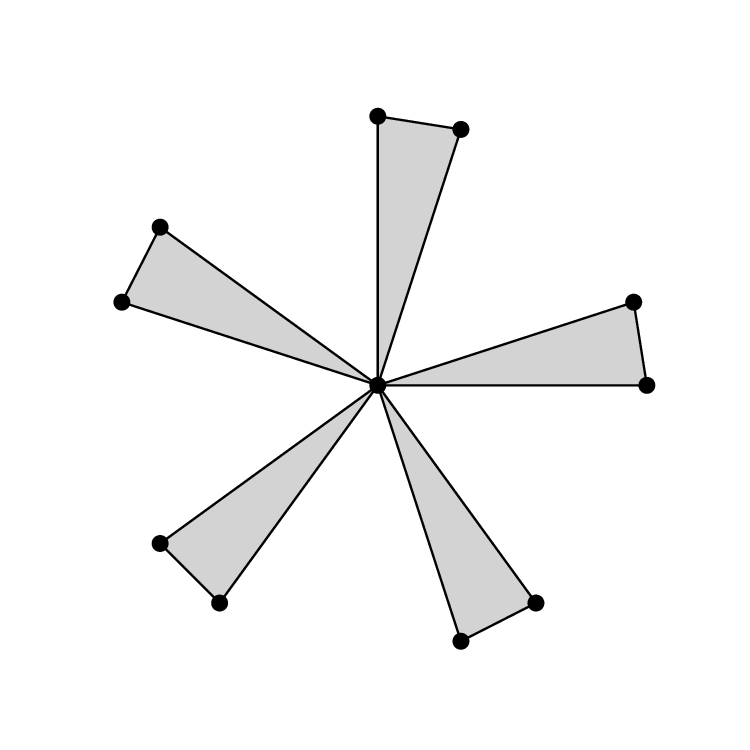}
			\caption{The $3$-uniform hyperstar with $5$ hyperedges}
			\label{fig2}
		\end{figure}
	\end{minipage}
\end{minipage}
\hspace{0.5cm}

Yue, Zhang and Lu \cite{yue2016largest} have obtained the explicit expression for the largest eigevalue of a $h$-uniform hyperpath of size $m=3,4$ and, for $m \geq 5$ they have given a polynomial whose largest real root will be the spectral radius of the hyperpath.

\begin{theorem}\cite{li2016extremal,zhang2018spectral}
    If  $\mathcal{T} \in \prescript{}{m}{T}^{(h)}$, then $2 \cos{\frac{\pi}{m+2}} \leq \rho(\mathcal{T}) \leq m^{\frac{1}{h}},$ equality is attained by the $h$-uniform hyperstar $\mathcal{S}_m^{(h)}.$
   Upper bound equality is reached if and only if $\mathcal{T} \cong \mathcal{S}_m^{(h)}$ and the attainment of the lower bound is possible if and only if $\mathcal{T} \cong \mathcal{P}_m^{(h)}.$
\end{theorem}

Zhang, Li and Guo \cite{zhang2020uniform} determined the connected $h$-uniform hypergraph that attains second minimum value of the largest eigenvalue in the class of all hypergraphs of uniformity $h$ and size $m$.
Suppose that the $h$-uniform hyperpath $\mathcal{P}_{m-1}^{(h)}$ is defined by the sequence $\mathcal{P}_{m-1}^{(h)}:=e_1v_1e_2v_2\ldots v_{m-2}e_{m-1},$ where $v_i \in e_i \cap e_{i+1}$, for $1 \leq i \leq m-2,$ and $\vert e_i \vert = r.$
We denote by  $\mathcal{D}_{m,h}$, the hypergraph obtained from $\mathcal{P}_{m-1}^{(h)}$ by the addition of a new pendent hyperedge $e_m$ to a pendent vertex of $e_2$.
Xiao and Wang \cite{xiao2021effect} have thoroughly examined the impact of the grafting operation \cite{li2016extremal} on the hypergraph's spectral radius.

\begin{theorem}\cite{zhang2020uniform,xiao2021effect,su2018matching}
    If $\mathcal{H} \in H^{(h)}$ is connected and $ \mathcal{H} \not\cong \mathcal{P}_{m}^{(h)}$. Then 
    $\rho({\mathcal{H}}) \leq \rho(\mathcal{D}_{m,h}).$
\end{theorem}

Wang and Yuan \cite{wang2020uniform} have determined the hypertree that attains third minimum value of the largest eigenvalue in the class of all uniform hypertrees of fixed size. 
The authors of the article \cite{yuan2016ordering} have defined new operation of \emph{edge-moving} and by using both edge-moving operations, they have ordered first eight $h$-uniform hypertrees (they call it supertrees) of given order having maximum spectral radius. Later, the hypertrees that attain ninth and tenth highest spectral radius are characterized by Yuan, Xi and Zhang \cite{yuan2017ordering}.

A vertex subset $I \subset \mathcal{V}$ of $\mathcal{H}$ is called an independent set (a stable set) if for any $e \in \mathcal{E}$, $e \not\subset I$.  
If $I$ denotes the maximum independent set in $\mathcal{H}$, then $\vert I \vert$ is the independence number (stability number) of $\mathcal{H}$.
Su and Li \cite{su2023hypertree} have determined an upper bound for the largest eigenvalue of a hypertree using the independence number and they have characterized the hypergraph that attains this bound.

\begin{definition}
    Let $\mathcal{S}_q^{(h)}$ be a $h$-uniform hyperstar with $q$ hyperedges, say $e_1,\ldots,e_q$. Now, $\mathcal{S}(b_1, \ldots, b_q), b_i \leq h-1 $ is a $h$-uniform hypergraph obtained from $\mathcal{S}_q^{(h)}$ by attaching $b_i$ disjoint hyperedges in the $b_i$ distinct core vertices of $e_i$, ~ $1 \leq i \leq q$. 
    If $b_1=\cdots=b_s, s < q $, then we write $\mathcal{S}(b_1, \ldots, b_q)$ as $\mathcal{S}(b_1^{(s)},b_{s+1}, \ldots, b_q).$
\end{definition}

\begin{theorem}\cite{su2023hypertree}
    Let $\mathcal{T}$ be a $h$-uniform hypertree with $m (> h)$ hyperedges having independence number $\alpha$. 
    \begin{equation*}
        \rho(\mathcal{T}) \leq \left( \frac{1}{1-\alpha*}\right)^{\frac{1}{h}}.
    \end{equation*}
    Here, $\alpha*$ is the maximum root of $y^{h-1}\left( \frac{1}{1-y} - \frac{1}{y^{t_3}} -t_2 \right) =t_1,$
    where $t_1,t_2,t_3$ are integers satisfying the equations 
    \begin{align*}
        n-1-\alpha  &= t_1(h-1) + t_3, \text{ where } 0 \leq t_3 < h-1, \\
        m           &= t_1r+t_2+t_3+1.
    \end{align*}
   Furthermore, equality holds in the upper bound if $\mathcal{T}$ is isomorphic to $\mathcal{S}((h-1)^{(t_1)}, t_3, 0^{(t_2)})).$
\end{theorem}

The hypertree attaining the maximum value of the largest eigenvalue in the class of $h$-uniform hypertrees of given order, size and independence number was also characterized in \cite{zhang2022extremal}.
In the following definition, the hypergraph that attaining the maximum value of the largest eigenvalue among hypertrees of fixed diameter is defined.

\begin{definition}
    Let $\mathcal{H} \in {H}^{(h)}$. For $v,u_1,u_2 \in \mathcal{V}(\mathcal{H})$, let  $\mathcal{H}_{t_1,t_2}(v) \in {H}^{(h)}$ be obtained from $\mathcal{H}$ by attaching (a core vertex of a pendent hyperedge of) hyperpaths $\mathcal{P}_{t_1}^{(h)}$ and $\mathcal{P}_{t_2}^{(h)}$ at $v$. Similarly, let  $\mathcal{H}_{t_1,t_2}(u_1,u_2) \in {H}^{(h)}$ be  obtained from $\mathcal{H}$ by attaching (a core vertex of a pendent hyperedge of) hyperpaths $\mathcal{P}_{t_1}^{(h)}$ and $\mathcal{P}_{t_2}^{(h)}$ at $u_1$ and $u_2,$ respectively.
\end{definition}

\begin{lemma}\cite{wang2020some}
    Given that $u_1,u_2$ are the non-pendent vertices of $\mathcal{H}$, and for any $t_1,t_2 \geq 1$,
    if the existence of the internal path $\mathcal{P}_s$ (of length $s$) in $\mathcal{H}_{t_1,t_2}(u_1,u_2)$  is guaranteed, then 
    \begin{equation*}
    {\rho(\mathcal{H}_{t_1+1,t_2-1}(u_1,u_2)) < \rho(\mathcal{H}_{t_1,t_2}(u_1,u_2)), \text{ for } t_1-t_2+1 \geq s \geq 0}. 
    \end{equation*}
\end{lemma}

\begin{lemma}\cite{guo2020alpha}
     Let $\mathcal{H} \in \mathcal{H}^{(r)}$ with at least one hyperedge and $v \in \mathcal{V}(\mathcal{H})$. For $t_1,t_2 \geq 1$, we have
    \begin{equation*}
    \rho(\mathcal{H}_{t_1+1,t_2-1}(v)) < \rho(\mathcal{H}_{t_1,t_2}(v)).
    \end{equation*}
\end{lemma}

\begin{definition}\cite{niu2021sharp}
    A hyperbug $\mathcal{B}_{t_1,t_2,t_3}$ is a $h$-uniform hypergraph obtained from an edge deleted  $h$-uniform complete hypergraph $\mathcal{K}_{t_1}^{(h)} \setminus e$, by attaching two hyperpaths $\mathcal{P}_{t_2}^{(h)}$ and $\mathcal{P}_{t_3}^{(h)}$ at $u_1$ and $u_2$, respectively, where $u_1,u_2 \in e.$ A hyperbug is said to be balanced if $\vert t_2 - t_3 \vert \leq 1.$
    \end{definition}

By using the above lemmas, Niu, Ren and Zhang could able to characterize the extremal hypergraphs attaining the maximum value of the largest eigenvalue when diameter is fixed, and the minimum value of the largest eigenvalue when clique number is fixed.

\begin{theorem}\cite{niu2021sharp}
 If $\mathcal{H} \in \prescript{n}{}{H}^{(h)}$ has diameter at least equal to $D$, then 
 \begin{equation*}
     \rho(\mathcal{H} ) \leq \rho \left(\mathcal{B}_{n-(D-2)(h-1), \left\lfloor \frac{D}{2} \right\rfloor, \left\lceil \frac{D}{2} \right\rceil})\right.
 \end{equation*}
The necessity and sufficiency for the equality is $\mathcal{H} \cong  \mathcal{B}_{n-(D-2)(h-1), \left\lfloor \frac{D}{2} \right\rfloor, \left\lceil \frac{D}{2} \right\rceil}$.
\end{theorem}

\begin{definition}\cite{niu2021sharp}
    A $h$-uniform kite hypergraph $\mathcal{KT}_{t_1,t_2}$ is obtained  from a $h$-uniform complete hypergraph $\mathcal{K}_{t_1}$ by attaching to any one of the vertex, a hyperpath $\mathcal{P}_{t_2}$.
\end{definition}

\begin{theorem}\cite{niu2021sharp}
   Let $\omega \geq 2$ be the clique number of $\mathcal{H} \in \prescript{n}{}{H}^{(h)}$. Then 
 \begin{equation*}
     \rho(\mathcal{H} ) \geq \rho( \mathcal{KT}_{\omega, n-\omega}).
 \end{equation*}
 The necessity and sufficiency for the equality is $\mathcal{H} \cong \mathcal{KT}_{\omega, n-\omega}.$
\end{theorem}

Xiao, Wang and Lu \cite{xiao2017maximum} have generalized the  BFS (\emph{Breadth-First-Search})  ordering of trees to uniform hypertrees, using which they characterized the extremal hypertree attaining the maximum value of the largest eigenvalue among the hypertrees of fixed degree sequence. The similar type of result in the case of trees was initially studied by B{\i}y{\i}koglu and Leydold \cite{biyikoglu2008graphs}, using the same technique. 

 Given the family of $h$-uniform hypergraphs with diameter $d$ and size $m$, the first two hypergraphs that attains the largest spectral radius was characterized by  Xiao, Wang and Du \cite{xiao2018first}. 

Let $\mathcal{P}_{m_j}^{(h)}$ be the $h$-uniform hyperpaths with $m_j$ hyperedges. 
The hyperpaths $\mathcal{P}_{m_1}^{(h)}, \ldots, \mathcal{P}_{m_p}^{(h)}$ are said to have almost equal lengths if the difference of their lengths is at most one.

\begin{theorem}\cite{xiao2019maximum,zhang2019maximum}
If  $\mathcal{T}^*(m,p,h) \in \prescript{}{m}{T}^{(h)}$ with $m=\sum\limits_{j=1}^p m_j$ denotes a hypertree consisting of $p$ pendent hyperpaths $\mathcal{P}_{m_j}^{(h)},1 \leq j \leq p$ of almost equal length, and these $p$  hyperpaths share  a common vertex. Then for any hypertree $\mathcal{T} \in \prescript{}{m}{T}^{(h)}$ with  $p$ pendent hyperedges 
$$ \rho(\mathcal{T}) \leq \rho(\mathcal{T}^*(m,p,h)),$$ 
and the necessity and sufficiency for the equality is $\mathcal{T} \cong \mathcal{T}^*(m,p,h)$.
\end{theorem}

\subsubsection{Matching}

\begin{definition}
    A matching $M$ in a hypergraph $\mathcal{H}$ is a set of independent hyperedges of $\mathcal{H}$, and the size of $M$ is the number of hyperedges in it ($\vert M \vert$).
    A $t$-matching of $\mathcal{H}$ is a matching of size $t.$  
    A maximal matching in $\mathcal{H}$ is a matching $M$ in $\mathcal{H}$ such that $M \cup \{ e \}$ is not a matching in $\mathcal{H}$ for any $e \in (\mathcal{E}\setminus M).$
    A maximum matching in $\mathcal{H}$ is a maximal matching of maximum cardinality. The matching number of $\mathcal{H}$, denoted by $\nu(\mathcal{H})$ is the size of a maximum matching.
\end{definition}

A detailed study of the hypergraph matching problem under more general conditions was carried out in \cite{keevash2015geometric} (but with slightly varying notations).
Using the BFS ordering of the uniform hypertrees, Guo and Zhou \cite{guo2018spectral} have determined the extremal hypertrees that attains the maximum value of the largest eigenvalue in the class of uniform hypertrees of given matching number \textbackslash  independence number \textbackslash number of branch vertices.

\begin{theorem}\cite{su2020largest}
    Let $\mathcal{T} \in \prescript{}{m}{T}^{(r)}$, $ m> h$ be having matching number $\nu$. 
    \begin{equation*}
        \rho(\mathcal{T}) \leq \left( \frac{1}{1-\alpha*}\right)^{\frac{1}{h}}.
    \end{equation*}
    Here, $\alpha*$ is the maximum root of $y^{h-1}\left( \frac{1}{1-y} - \frac{1}{y^{t_3}} -t_2 \right) =t_1,$
    where $t_1,t_2$ and  $t_3$ are integers satisfying the equations 
    \begin{align*}
        \nu - 1 &= (h-1)t_1 +t_3, \text{ where } 0 \leq t_3 < h-1, \\ 
              m &= t_1r+t_2+t_3+1.
    \end{align*}
    Furthermore, the attainment of the equality is when \\ $\mathcal{T} \cong \mathcal{S}((h-1)^{(t_1)}, t_3, 0^{(t_2)})).$
\end{theorem}

The hypertree attaining second maximum value of the largest eigenvalue in the class of hypertrees of fixed size and matching number was characterized by Su et al. \cite{su2021second}. In the same year,
Yu at al. \cite{yu2021spectral} have completely characterized the hypergraph attaining the maximum value of the largest eigenvalue in the class of unicyclic uniform hypergraphs of fixed matching number.

\noindent A matching is said to be a perfect matching if it covers all the vertices of the hypergraph. That is, $\bigcup\limits_{e \in M} e = \mathcal{V}(\mathcal{H}).$
Zhang and Chang \cite{zhang2018spectral} have given the upper bound for the spectral radius of the $h$-uniform hyprtrees with perfect matching and they have also characterized the extremal graph achieving this upper bound. 

\begin{definition}\label{extremaltree_perfectmatching_spectralradius}
    Let $\mathcal{T}_1 \in T^{(h)}$ be obtained from a hyperstar (on $n_1=m_1(h-1)+1$ vertices) $\mathcal{S}_{m_1}^{(h)}$ by attaching a pendent hyperedge at each vertex of $\mathcal{S}_{m_1}^{(h)}.$
\end{definition}

\begin{theorem}\cite{xu1997spectral}
    Let $\mathcal{T} \in \prescript{2t}{}{T}^{(2)},t\geq 1$. Then 
    \begin{equation*}
        \rho(\mathcal{T}) \leq \frac{1}{2}(\sqrt{t-1}+\sqrt{t-3}),
    \end{equation*}
    and the equality holds if and only if $T$ is isomorphic to $\mathcal{T}_1$ (which is a $2$-uniform hypertree defined in Definition \ref{extremaltree_perfectmatching_spectralradius}).
\end{theorem}

\begin{theorem}\cite{zhang2018spectral}
    Let $\mathcal{T} \in \prescript{ht}{}{T}^{(h)}$. Then, $\rho(\mathcal{T}) \leq \rho(\mathcal{T}_1)$ ($\mathcal{T}_1$ is a $h$-uniform hypertree defined in Definition \ref{extremaltree_perfectmatching_spectralradius}) and 
    $\rho(\mathcal{T}_1)$ is a maximum real root of the equation (in the variable $x$)
    \begin{equation*}
        x^h-(t-1)^{\frac{1}{h}}x^{h-1}-1=0.
    \end{equation*}
\end{theorem}

Among all $h$-uniform unicyclic (linear and non-linear case) hypergraphs having perfect matching, Sun, Wang, and Ni \cite{suna2023maximal} have characterized the extremal hypergraph that attains the maximum spectral radius.

\section{Characteristic polynomial}

In literature, the computation of the characteristic polynomial of the various hypermatrices associated with the hypergraph was carried out through either of the techniques.
The first one is through Poisson's formula \cite{cox2005using}; while the other is through the generalized traces \cite{morozov2011analogue}.

 Chen and Bu \cite{chen2021reduction} have obtained the multiplicity of the point zero ($R_1=(h-1)^{h-1}-h^{h-2}$) and the sum of the multiplicities of the non-zero points ($R_2=h^{h-2}$) in the affine variety (defined by the $f_i's$), with the help of the existing expression for the characteristic polynomial of a single hyperedge \cite{cooper2012spectra}.

\begin{theorem}\cite{cooper2012spectra}
    If $\mathcal{P}_1^{(h)}$ denotes a hypergraph with single hyperedge of cardinality $h$, then
    $$\Phi_{P_1^{(h)}}(\lambda) =  \lambda^{h(h-1)^{h-1}-h^{h-1}}(\lambda^h-1)^{h^{h-2}}.$$
\end{theorem}

\begin{theorem}\cite{cooper2015computing}
    If $\mathcal{S}^{(3)}_{m}$ denotes a hyperstar with uniformity $3$ and size $m$, then 
    \begin{equation*}
      \Phi_{S^{(3)}_{m}}(\lambda) =  \lambda^{(2m-2)4^{m}} \prod\limits_{t=0}^m (\lambda^3-t)^{\binom{m}{t}3^t}.
    \end{equation*}
\end{theorem}

\begin{theorem}\cite{bao2020combinatorial}
    If $S^{(h)}_{m}$ denotes a hyperstar with uniformity $h$ and size $m$, then 
    $$ \Phi_{S^{(h)}_{m}}(\lambda) = \lambda^{m(h-1)^{m(h-1)+1}} \prod\limits_{i=1}^m \left( \lambda -\frac{i}{\lambda^{h-1}} \right)^{\binom{m}{i}h^{(h-2)i}((h-1)^{h-1}-h^{h-2})^{m-i}}.$$
\end{theorem}

\begin{theorem}\cite{duan2023characteristic}
    Let $S^{(h)}_{m_1,m_2}$ be a $h$-uniform double hyperstar with $m_1+m_2+1$ hyperedges. Then,
    $ \Phi_{S^{(h)}_{m_1,m_2}}(\lambda) = \lambda^{((m-2)(h-1)-1)(h-1)^{m(h-1)}} \prod\limits_{p=0}^{m_1} \prod\limits_{q=0}^{m_2} \\(\lambda^h-p)^{\binom{m_1}{p} R_1^p R_2^{m_1-p+1} (h-1)^{m_2(h-1)}}(\lambda^h-q)^{\binom{m_2}{q} R_1^q R_2^{m_2-q+1} (h-1)^{m_1(h-1)} } (\lambda^{2r}-  (p+q+1)\lambda^h+pq)^{\binom{m_2}{q}R_1^{q+1}R_2^{m_2-q} \binom{m_1}{p}R_1^p R_2^{m_1-p}},$ 
    where $R_1=h^{h-2}$ and $R_2=(h-1)^{h-1}-h^{h-2}.$   
\end{theorem}

\begin{theorem}\cite{chen2021reduction}
    If  $P_m^{(h)}$ denotes hyperpath with uniformity $h$ and size $m$, then
    $$ \Phi_{P_m^{(h)}}(\lambda) = \prod\limits_{i=0}^m \Phi_{P_i^{(2)}}(\lambda^{\frac{h}{2}})^{c(i,m)}, \text{ where }$$
     $c(i,m)=\begin{cases}
        R_2^m, &  j=m \\
        ((m-i-1)R_1+2R_2)R_1R_2^i(R_1+R_2)^{m-i-2}, &   1\leq i \leq m-1 \\
        \frac{2}{h}[m(h-1)+1](R_1+R_2)^m, & i=0
    \end{cases},$ \\ $ R_1=(h-1)^{h-1}-h^{h-2}$ and $R_2=h^{h-2}.$
\end{theorem}

\begin{theorem}\cite{chen2021reduction}
Suppose that $\mathcal{H} \in \prescript{n}{}{H}^{(h)}$.
Let $\mathcal{H}_u^s \in {H}^{(h)}$ be obtained from $\mathcal{H}$ by adding $s$ pendent hyperedges at the vertex $u$ of $\mathcal{H}$. Then 
    $$ \Phi_{\mathcal{H}_u^s}(\lambda) = \lambda^{s(h-1)^{n+s(h-1)}} \Phi_{\mathcal{H}-v}(\lambda)^{(h-1)^{s(h-1)+1}} \prod\limits_{p=0}^s \left( M_{\mathcal{H}}\left( \lambda, \frac{1}{\lambda} \right) \right)^{\binom{s}{p}R_1^{s-p}R_2^p},$$
    where $R_1=(h-1)^{h-1}-h^{h-2}$ and $R_2=h^{h-2}.$
\end{theorem}

Cooper and Dutle \cite{cooper2012spectra} found that the $0,1$ and $\binom{n-1}{2}$ are the eigenvalues of a $3$-uniform complete hypergraph.
 Zheng \cite{zheng2020simplifying} proposed a method that simplifies the computation of the eigenvalues of $h$-uniform ($h=3,4$) complete hypergraph of given order.
 In \cite{zheng2021characteristic}, he has given the explicit expression for the characteristic polynomial of complete $3$-uniform hypergraph on $n$ vertices using the Poisson formula for resultants and hence the algebraic multiplicity of the eigenvalues $0,1$ and $\binom{n-1}{2}$ is determined.

\begin{theorem}\cite{zheng2021characteristic}
    Let $\mathcal{H} \in \prescript{n}{}{H}^{(h)},n \geq 3$ be a complete hypergraph. Then\\
    1. $\Phi(\mathcal{H})=    (\lambda-1)^{(n-3)2^{n-1}} \prod\limits_{t=0}^{n-1} \left( \lambda^3 - \frac{1}{2}[(n-2t)^2 -3(n-2)]\lambda^2 \right.\\ \left. + \frac{1}{2}[(n-2)^2 + (n-t-1)^2+(t-1)^2]\lambda + \frac{1}{2}(n-2)(n-t-1)(t-1) \right)^{\binom{n-1}{t}}$. \\
        2. $am_{\mathcal{H}}(0)=n;am_{\mathcal{H}}(\binom{n-1}{2})=1;am_{\mathcal{H}}(1)=(n-2)2^{n-1}+2.$
\end{theorem}

\begin{theorem}\cite{duan2023characteristic}
Let $\mathcal{C}_3^{(h)}$ denote the $h$-uniform hypercycle of length $3$. Then,
$$\Phi_{\mathcal{C}_3^{(h)}}(\lambda) = \lambda^{c_o} (\lambda^h-1)^{c_1} (\lambda^h-2)^{c_2}(\lambda^h-4)^{c_4},$$
where $c_0=3(h-1)^{3h-3}-3h^{h-1}(h-1)^{2h-3}+3h^{2h-3}(h-1)^{h-2}-3h^{3h-6},\\ 
c_1 = 3h^{h-2}(h-1)^{2h-3}-6h^{2h-4}(h-1)^{h-2}+8h^{3h-7}, \\
c_2= 3h^{2h-4}(h-1)^{h-2}-6h^{3h-7}, c_4=h^{3h-7}$
\end{theorem}

\begin{theorem}\cite{duan2024characteristic}
Let $\mathcal{C}_4^{(h)}$ denote the $h$-uniform hypercycle of length $4$. Then, $\Phi_{\mathcal{C}_4^{(h)}}(\lambda) =$
$$\lambda^{c_o} (\lambda^h-1)^{c_1} (\lambda^h-2)^{c_2}(\lambda^h-4)^{c_4}\left(\lambda^h-\frac{3+\sqrt{5}}{2}\right)^{c'}\left(\lambda^h-\frac{3+\sqrt{5}}{2}\right)^{c'},$$
where $c_0=4(h-1)^{4h-4}-4h^{h-1}(h-1)^{3h-4}+4h^{2h-3}(h-1)^{2h-3}-4h^{3h-5}(h-1)^{h-2}+5h^{4h-8}, \\
c_1 =4h^{h-2}(h-1)^{3h-4}-8h^{2h-4}(h-1)^{2h-3}+4h^{3h-6}(h-1)^{h-2},\\
c_2=4h^{2h-4}(h-1)^{2h-3}-8h^{3h-6}(h-1)^{h-2}+ 10h^{4h-9}, \\
c_4=h^{4h-9}, c'=4h^{3h-5}(h-1)^{h-2}-8h^{4h-9}. $
\end{theorem}

Let $\mathcal{B}=\{ \beta_i^2 : \beta_i \text{ is an eigenvalue of the signed subgraph of } G \}$ be the set of squares of all (non-zero) eigenvalues of all the signed subgraphs of $G$. Also, let $\vert \mathcal{B} \vert= s$ and $b_1,\ldots,b_s$ be the distinct elements of $\mathcal{B}.$
Let $M$ be an $s\times s$ matrix, whose $ij^{th}$ entry ($1 \leq i,j \leq s$) will be equal to $(b_j)^i$ (that is, the $i^{th}$ power of $b_j$). 
Since $b_i$'s are all distinct and we know that the (square) Vandermonde matrix is invertible if and only if all the entries in the first row are distinct. 
Let $\mathcal{G}(s)$ denote the set of all non-isomorphic subgraphs of $G$ with at most $s$ edges. Suppose that $\vert \mathcal{G}(s) \vert = k$ and $G_1,\ldots, G_k$ are the elements of $\mathcal{G}(s).$
A closed walk that uses every edge an even number of times is called a parity-closed walk and  a parity-closed walk of length $i$ in $G$ is denoted by $p_i(G)$.
Let $P$ be a $k \times k$ (parity-closed walk) matrix, that is the $ij^{th}$ entry of the matrix $P$ is $p_{2i}(G_j)$. 
Let $N=(N_G(G_1),\ldots,N_G(G_s))^\intercal$ be a column vector length $s$, where $N_G(G_i)$ denotes the number of subgraphs $G'$ of $G$ that are isomorphic to $G_i.$
Also, let $D(h)$ be the diagonal matrix with 
$(D(h))_{ii} = \frac{2^{m_i-n_i}k^{m_i(k-3) + n_i}}{(k-1)^{n_i+m_i(k-2)-1}}, $ where $m_i$ and $n_i$ are respectively, the number of edges and vertices in the subgraph $G_i$.

\begin{theorem}\cite{chen2024spectra}
    Given  $h \geq 3$ and a graph of order $n'$ and size $m$. The characteristic polynomial of the power hypergraph $G^{h}$ is given by 
    \begin{equation*}
        \Phi_{G^{h}}(\lambda) = \lambda^{c_0(h)}\prod\limits_{i=1}^{s}(\lambda^k-b_i)^{c_i(h)},
    \end{equation*}
    where $c_i(h) = \frac{(h-1)^{n+(h-2)m-1}}{h}(M^{-1}PD(h)N)_i$ for $1 \leq i \leq s$, and $c_0=(n'+(h-2)m)(h-1)^{n'+(h-2)m-1}-r\sum\limits_{i=1}^s c_i(h)$ .
\end{theorem}

A subgraph $G_1$ of a graph $G$ is said to be an elementary subgraph if every component (maximal connected subgraph) of $G_1$ is either a cycle or an edge.
Initially, Harary \cite{harary1962determinant} in the year 1962 has shown that the determinant of adjacency matrix of the graph can be obtained from the number of elementary subgraphs of the  graph. Later, Sachs \cite{{sachs1966teiler}} in 1966 has given more explicit formula for the characteristic polynomial of the graph in terms of the number of elementary subgraphs.

\begin{theorem}\cite{harary1962determinant,sachs1966teiler}[Harary-Sachs Theorem for Graphs]
Given a simple labeled graph $G$ on $n$ vertices. Let $c_j$ be the co-degree $j$ (i.e., $\lambda^{n-j}$) co-efficient of the characteristic polynomial of $G$ 
and let $G_i$ be the set of all subgraphs of $G$ on $i$ vertices with all its components are either an edge or a cycle. 
Then, 
$$c_j = \sum\limits_{G' \in G_j} (-1)^{c_1(G')+c_2(G')}2^{c_2(G')} [\#G'\subseteq G],$$
where $c_1(G')$ and $c_2(G')$ are respectively, the number of edge and cycle components in $G'$, and  $[\#G'\subseteq G]$ denotes the number of labeled subgraphs of $G$ that are isomorphic to $G'$.
\end{theorem}

\noindent The generalization of the Harary-Sachs theorem from graphs to a $h$-uniform hypergraphs can be found in \cite{clark2021harary}. By keeping the similar notations, the result is formally stated in the following theorem.

\begin{definition}
    The \emph{flattening} of a $h$-uniform, labeled, multi-hypergraph $\mathbf{H}$ is the removal of the duplicate hyperedges, (that is, the multiplicity of every hyperedge is reduced to one) and the resulting simple hypergraph is denoted by $\underline{\mathbf{H}}.$
    We call $\mathbf{H}$ is an infragraph of $\mathcal{H}$ if $\underline{\mathbf{H}}$ is a sub-hypergraph of $\mathcal{H}.$
\end{definition}

\begin{definition}
\begin{itemize}
    \item A $u$-rooted star of a hyperedge $e$, denoted by $S_{e}(u)$ is a digraph with vertex set $e$ and the arc set is given by $\{ (u,v): v \in e, v \neq u \}.$
    \item An ordering $R=(S_{e_1}(u_1), \ldots, S_{e_m}(u_m))$ of the rooted star such that $\mathcal{E}(\mathcal{H})=\{e_1,\ldots,e_m \}$ and $u_i \leq u_{i=1},$ is called a rooting of a $h$-uniform hypergraph $\mathcal{H}$. The rooted multi-digraph corresponding to a rooting $R$ of $\mathcal{H},$ denoted by $D(R)$ or simply $D$ is defined as $D = \bigcup\limits_{i=1}^{m} S_{e_m}(u_m),$ where the union adds up the multiplicity of the arcs.
    \item If a rooted multi-digraph $D$ corresponding to a rooting $R$ of a hypergraph $\mathcal{H}$ is called an Euler rooting if $D$ is Eulerian. We denote by $\mathcal{R}(\mathcal{H}),$ the multi-set of all Euler rooted digraphs of $\mathcal{H}.$
    \end{itemize}
\end{definition}

\begin{figure}[h]
    \centering
    \includegraphics[width=0.4\linewidth]{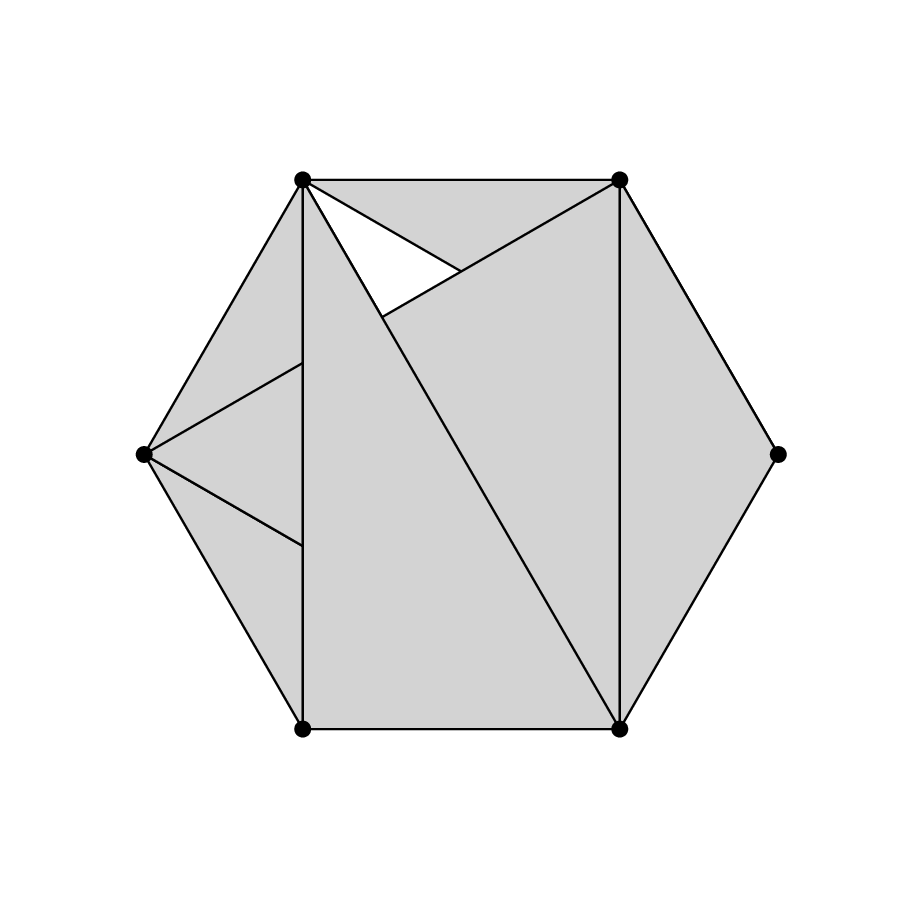}
    \caption{$3$-uniform Veblen hypergraph}
    \label{veblen_hypergraph}
\end{figure}

\begin{definition}
   A hypergraph is said to be $h$-valent if the degree of every vertex is a multiple of $h$.
   A $h$-uniform, $h$-valent (multi) hypergraph is the \emph{Veblen} hypergraph.
\end{definition}
\noindent A $3$-uniform Veblen hypergraph with $6$ hyperedges is pictorially depicted in Figure \ref{veblen_hypergraph}.

\begin{theorem}\cite{veblen1912application}[Veblen's Theorem]
    A (multi) graph can be decomposed into disjoint cycles if and only if the degree of every vertex is even.
\end{theorem}

\begin{theorem}\cite{van1987circuits}[BEST Theorem]
 The number of Euler circuits in a connected Eulerian graph $G$ is 
 \begin{equation*}
     |E(G)|=\tau(G)\prod\limits_{u \in \mathcal{V}(G)} (d_v-1)!
 \end{equation*} 
\end{theorem}

\begin{theorem}\cite{clark2021harary}[Harary-Sachs Theorem for Hypergraphs]
 Given a simple $h$-uniform hypergraph $\mathcal{H}$ on $n$ vertices. Let $\nu_d^*(\mathcal{H})$ denotes the set of Veblen infragraphs of $\mathcal{H}$ with $d$ hyperedges up to isomorphism.
 If $c_d$ denotes the co-degree $d$ co-efficient of the characteristic polynomial of $\mathcal{H},$ then
\begin{equation}
    c_d= \sum\limits_{\mathbf{H} \in \nu_d^*(\mathcal{H}) } (-1)^{c(\mathbf{H})} \left((h-1)^{n} \right)^{c(\mathbf{H})} \mathcal{C}_{\mathbf{H}}(\# \mathbf{H} \subseteq \mathcal{H}),
\end{equation}
where
$c(\mathbf{H})$ is the number of components in $\mathbf{H},$ and if $\mathbf{H} = \bigcup\limits_{i=1}^{c(\mathbf{H})} \Gamma_i,$ then 
$$(\# \mathbf{H} \subseteq \mathcal{H})= \frac{\mu_{\mathbf{H}}}{{c(\mathbf{H})}!} \prod\limits_{i=1}^{c(\mathbf{H})} ( \# \Gamma_i \subseteq \mathcal{H}) \text{ and } \mathcal{C}_{\mathbf{H}} = \prod\limits_{i=1}^{c(\mathbf{H})} \mathcal{C}_{\Gamma_i},  $$ \text{ where } $$( \# \Gamma_i \subseteq \mathcal{H}) = \left\vert {Aut(\underline{\Gamma_i})}  / {Aut({\Gamma_i})} \right\vert | \{ S \subseteq \mathcal{H}: S \cong \underline{\Gamma_i}  \} | $$ 
$$\text{ and } \mathcal{C}_{\Gamma_i}= \sum\limits_{ D \in \mathcal{R}(\Gamma_i)} \frac{\tau(D)}{\prod\limits_{u \in V(D)} d^-(u)}. $$ \\
Here $ \mu_{\mathbf{H}}= \binom{c(\mathbf{H})}{n_1,\ldots,n_t},$ where $t$ is the number of isomorphism classes (that is, the components of $\mathbf{H}$ that are isomorphic are belong to the same class) and $n_j$ be the cardinality of those classes, $d^-(u)$ denotes the in degree of a vertex $u$ in the digraph $D$ and  $\tau(D)$ denotes the number of arborescence in  $D,$ where $\mathcal{R}(\mathcal{H})$ is the multi-set of all the Euler rooted digraphs of $\mathcal{H}.$
\end{theorem}

\noindent Clark and Cooper  \cite{clark2022applications} have also shown that the above generalization of the Harary-Sachs theorem for hypergraphs is a faithful generalization (by showing that the result reduces to the classical Harary-Sachs theorem if we take $h=2$ in the above theorem) and also discussed some applications of the above stated theorem.

\subsection{Polynomial reconstruction problem}

Graph reconstruction is a problem of constructing the original graph $G$, where the multi-set of all vertex deleted subgraphs (also called vertex deck) of $G$ or edge deleted subgraph of $G$.
The conjecture \cite{ulam1960collection,kelly1957congruence} was that any graph of order at least $3$ is uniquely determined (up to isomorphism) by its vertex deck. 
In the case of trees, the conjecture was proved by Kelly \cite{kelly1957congruence}, and for any triangle free graphs (on putting some constraints on diameter and connectivity) by Clifton et al. \cite{clifton2024reconstruction}. The above conjecture is still open for general simple graphs. 
Later, the polynomial version of the conjecture was posed, and the problem of polynomial reconstruction was studied in \cite{schwenk1979spectral,cvetkovic1988recent}.
Polynomial reconstruction of a graph is a problem of constructing the characteristic polynomial of the graph $G$, where the characteristic polynomial of all the vertex deleted subgraphs of $G$ are known. Some of the recent developments in this direction can be found in \cite{farrugia2021graphs,sciriha2023polynomial}. 
It was proved \cite{tutte1979all} that, from the multi-set of all the vertex-deleted subgraphs of $G$ the characteristic polynomial $G$ can be constructed. 
It was also shown \cite{sciriha2023polynomial} that there exist non-isomorphic graphs whose multi-set of characteristic polynomials of the vertex-deleted subgraphs are the same.
Two hypergraphs $\mathcal{H}_1$ and $\mathcal{H}_2$ are said to be hypomorphic if there is a bijection $\theta'$ from $\mathcal{V}_1$ to $\mathcal{V}_2$ such that for any vertex $u \in \mathcal{V}_1$, $\mathcal{H}_1 - u \cong \mathcal{H}_2 - \theta'(u)$. 

The hypergraph reconstruction conjecture for $3$-uniform hypergraphs was disproved by Kocay \cite{kocay1987family} in 1987. 
He constructed infinite pairs $(\mathcal{G}_{j},\mathcal{H}_{j})_{j \geq 3}$ of non-isomorphic hypergraphs that are hypomorphic. 
Very recently, Cooper and Okur \cite{cooper2024polynomial} disproved the polynomial reconstruction conjecture for $3$-uniform hypergraphs by showing that the above pairs of hypergraphs have distinct characteristic polynomial.

\subsection{Relationship with the matching polynomial}

\begin{definition}\cite{zhang2017spectra}
An eigenvalue $\lambda$ of a dimension $n$ hypermatrix is said to be a \emph{totally non-zero eigenvalue} if the eigenvector corresponding to $\lambda$ has the support equal to $[n]$.
\end{definition}

\begin{theorem}\cite{zhang2017spectra}
    Let $M(\mathcal{T},i)$ denote the $i$-matching of a $h$-uniform hypertree $\mathcal{T}$ with matching number $\nu$, and 
    \begin{equation*}
        \phi_{\mathcal{T}}(x) = \sum\limits_{j=0}^m (-1)^j \vert M(\mathcal{T},j) \vert x^{(\nu-j)h}
    \end{equation*} be a polynomial in the indeterminate $x$. Then, $\lambda$ is a totally non-zero eigenvalue of $\mathcal{T}$ if and only if $\lambda$ is a root of the polynomial $\phi(\mathcal{T})$.
\end{theorem}

\begin{theorem}\cite{clark2018adjacency}
    Let $\mathcal{T}$  be a $h$-uniform hypertree. Then $\lambda$ is a nonzero eigenvalue of  $\mathcal{T}$ if and only if there exists an induced sub-hypertree $\mathcal{T}_1$ of $\mathcal{T}$ such that  $\lambda$ is a root of the polynomial  $\phi_{\mathcal{T}_1}(x).$
\end{theorem}

\begin{definition}\cite{su2018matching}
   A polynomial $\varphi_{\mathcal{H}}(x) = \sum\limits_{i \geq 0} (-1)^i a_{hi} x^{n-hi}$  associated with a $h$-uniform hypergraph $\mathcal{H}$ in the indeterminate $x$, is called a matching polynomial of $\mathcal{H}$ if $a_{2i}$ = is the number of $i$-matching of $\mathcal{H}$. 
\end{definition}

\begin{theorem}\cite{su2018matching}
   $1)~$ For any two vertex disjoint hypergraphs $\mathcal{H}_1$ and $ \mathcal{H}_2$, 
    \begin{align*}
          \varphi_{\mathcal{H}_1 \cup \mathcal{H}_2}(\lambda) =  \varphi_{\mathcal{H}_1}(\lambda) \varphi_{\mathcal{H}_2}(\lambda).
    \end{align*}
   $2)~ $ For a hyperedge $e$ of $\mathcal{H}$, if $\mathcal{H}- \mathcal{V}(e)$ denotes the hypergraph obtained from $\mathcal{H}$ by deleting the vertices of $e$, then
   \begin{equation*}
       \varphi_{\mathcal{H}}(\lambda) = \varphi_{\mathcal{H}-e}(\lambda) - \varphi_{\mathcal{H}-\mathcal{V}(e)}(\lambda).
   \end{equation*}
   $3)~$ Given a vertex $v$ of the hypergraph $\mathcal{H}$ and $I_1= \{ i : e_i \in \mathcal{E}_v(\mathcal{H}) \},$ then for any $I_2 \subseteq I_1$, we have 
   \begin{align*}
       \varphi_{\mathcal{H}}(\lambda) &= \varphi_{\mathcal{H}-\{ e_i : i \in I_2 \}}(\lambda) - \sum\limits_{i \in I_2} \varphi_{\mathcal{H} - \mathcal{V}(e_i)}(\lambda);\\
        \varphi_{\mathcal{H}}(\lambda) &= \lambda \varphi_{\mathcal{H}-u}(\lambda) - \sum\limits_{e \in \mathcal{E}_u(\mathcal{H})} \varphi_{\mathcal{H} - \mathcal{V}(e)}(\lambda).
   \end{align*}
\end{theorem}

\begin{theorem}\cite{wan2022distribution}
    The multi-set of all roots of the matching polynomial $\phi_{\mathcal{H}}(x)$ of a $h$-uniform hypergraph $\mathcal{H}$ is $h'$-symmetric for some positive integer $h'\leq h$.
\end{theorem}

\begin{theorem}\cite{su2018matching}
    Let $\mathcal{T} \in \prescript{}{m}{T}^{(h)}$.
    Then $\rho(\mathcal{T})$ is a simple root of the matching polynomial $\phi_{\mathcal{T}}(x)$ of $\mathcal{T}$.
\end{theorem}

\begin{theorem}\cite{mowshowitz1972characteristic}
    Let $\Phi_T(\lambda)$ and $\varphi_T(\lambda)$ respectively be the characteristic polynomial and matching polynomial of a tree $T$ ($2$-uniform hypertree). Then $\Phi_T(\lambda) = \varphi_T(\lambda).$
\end{theorem}

\begin{definition}
  A hyperedge $e$ of a connected sub-hypergraph $\mathbf{H}$ of $\mathcal{H}$ is said to be a boundary hyperedge if $e \cap \mathcal{V}(\mathbf{H}) \cap (\mathcal{V}(\mathcal{H}) \setminus \mathcal{V}(\mathbf{H})) \neq \emptyset .$ The set of all boundary hyperedges of the connected sub-hypergraph $\mathbf{T}$ is the boundary of $\mathbf{T}$ and is denoted by $\partial(\mathbf{T}).$  
\end{definition}

\begin{theorem}\cite{li2024relationship}
Let $\mathcal{T} \in \prescript{}{m}{T}^{(h)}$. Then 
\begin{equation*}
    \Phi_{\mathcal{T}}(\lambda)=\prod\limits_{\mathbf{T} \sqsubseteq \mathcal{T}} \varphi_{\mathbf{T}}(\lambda)^{a(\mathbf{T})},
\end{equation*}
where the product runs over all the connected sub-hypergraphs $\mathbf{T}$ of $\mathcal{T},$ and
\begin{equation*}
    a(\mathbf{T})= b^{m-e(\mathbf{T})-\vert \partial(\mathbf{T})  \vert} c^{e(\mathbf{T})} (b-c)^{\vert \partial(\mathbf{T})  \vert},
\end{equation*}
where $e(\mathbf{T})$ is the size of $\mathbf{T}$, $\partial(\mathbf{T})$ is the boundary of $\mathbf{T},$ $b=  (h-1)^{h-1} ,\ c=h^{h-2}$.
\end{theorem}

\subsection{Spectral moments and Estrada index of the hypergraph}

The $t^{th}$ order spectral moment of a hypermatrix $M$, denoted by $S_t(\mathcal{M})$ is defined to be the sum of $t^{th}$ power of all the eigenvalues of $\mathcal{M}.$ 
Hu et al. \cite{hu2013determinants} have shown that the  order $t$ trace of a hypermatrix $\mathcal{M}$ is equal to the $t^{th}$ spectral moment of $\mathcal{M},$ i.e., $Tr_{t}(\mathcal{M})=S_t(\mathcal{M}).$
The $t^{th}$ order spectral moment of a $h$-uniform hypergraph is the spectral moment of its adjacency hypermatrix.
The following characterization for  the symmetric spectrum of the hypergraph is by Shao, Qi and Hu.
\begin{theorem}\cite{shao2015some}
    Let $\Phi_{\mathcal{H}}(\lambda) = \sum\limits_{j=0}^d a_j \lambda^{d-j} $ be the characteristic polynomial of $\mathcal{H}$. Then the following are equivalent:
    \begin{itemize}
        \item The $\sigma(\mathcal{H})$ is $h$-symmetric;
        \item If $h\not| ~t$, then the coefficient $a_t=0;$
        \item If $h\not| ~t$, then $S_t(\mathcal{H})=0.$
    \end{itemize}
\end{theorem}

In the case graphs, it has been shown \cite{van2003graphs} that the characteristic polynomial of two graphs are same if and only if the order $t$ spectral moments of the two graphs are same for all $t$.

\begin{theorem} \cite{boulet2008lollipop,cvetkovic1987spectra,zhang2009lollipop}
For a connected subgraph $G'$ of a graph ${G},$  let $[\#G'\subseteq G]$ be the number of labeled subgraphs of $G$ that are isomorphic to $G'$,  and $c_t(G')$ be the number of closed walks of $G'$ with length $t$ in $G'$. Then 
\begin{equation*}
    S_t(G)= \sum\limits_{G' \in \mathcal{G}(t)} c_t(G') [\#G'\subseteq G],
\end{equation*}
where $\mathcal{G}(t)$ is the set of connected subgraphs of $G$ with at most $t$ edges.
\end{theorem}

Suppose that $h \vert t. $ 
Let $\mathcal{T}(\frac{t}{h},h)$ denote the set of $h$-uniform sub-hypertrees with at most $\frac{t}{h}$ hyperedges.
For a hypertree $\mathcal{T}' \in \mathcal{T}(\frac{t}{h}, h)$,  let $[\#\mathcal{T}'\subseteq \mathcal{T}]$ denote the number of sub-hypertrees of  $\mathcal{T}$ that are isomorphic to $\mathcal{T}'$.
    Here, the $t^{th}$ order spectral moment coefficient of the sub-hypertree $\mathcal{T}'$ is $$c_t(\mathcal{T}') = t \sum\limits_{\sum\limits_{e \in \mathcal{E}(\mathcal{T}')} w(e) = \frac{t}{h}} \left( \prod\limits_{e \in \mathcal{E}(\mathcal{T}')} w(e)^{h-1} \prod\limits_{v \in \mathcal{V}(\mathcal{T}')} \frac{(D_v-1)!}{R_v} \right),$$ where $D_v= \sum\limits_{e \in \mathcal{E}_v(\mathcal{T}')} w(e)$ and $R_v = \prod\limits_{ e \in \mathcal{E}_v(\mathcal{T}')} w(e)! $.

\begin{theorem} \cite{chen2022spectral}
    Let $\mathcal{T}$ be a $h$-uniform hypertree of size $m$ and $\mathcal{T}'$ be a sub-hypertree of $\mathcal{T}$ of size $m'$.  Then the $t^{th}$ spectral moment of $\mathcal{T}$ is,   $S_t(\mathcal{T})$ =  
    \begin{equation*}
      \begin{cases}
            {\sum\limits_{\mathcal{T}' \in \mathcal{T}(\frac{t}{h}, h)}   (h-1)^{(m-m')(h-1)} h^{m'(h-2)} c_t(\mathcal{T}')  [\#\mathcal{T}'\subseteq \mathcal{T}]} & \text{ if } h \vert t, \\
            0 & \text{ otherwise}.
        \end{cases}
    \end{equation*}  
\end{theorem}

In the following theorem the expression for the first $3h$ co-degree coefficients of the is obtained using the above result. 
\begin{theorem}\cite{chen2022spectral}
Let $\mathcal{P}_m^{(h)}$ and $\mathcal{S}_m^{(h)}$ be the $h$-uniform hyperpath and hyperstar, respectively.
If  $\Phi_{\mathcal{T}}(\lambda) = \sum\limits_{j=0}^d a_j \lambda^{d-j}$ denotes the characteristic polynomial of the hypertree $\mathcal{T} \in \prescript{n}{}{T}^{(h)}$, then for $0 \leq i \leq 3h$
\begin{equation*}
    a_i=\begin{cases}
        -h^{h-2} (h-1)^{n-r} [\#\mathcal{P}_1^{h}\subseteq \mathcal{T}] & \text{ if } i=h, \\
        \frac{1}{2} \left( \substack{h^{2h-4}(h-1)^{2(n-r)} [\#\mathcal{P}_1^{h}\subseteq \mathcal{T}]^2\\ -h^{h-2}(h-1)^{n-r}[\#\mathcal{P}_1^{h}\subseteq \mathcal{T}] \\- 2r{2h-4} (h-1)^{n-2r+1} [\#\mathcal{P}_2^{h}\subseteq \mathcal{T}]} \right) & \text{ if } i=2h, \\
        \frac{1}{6} \left( \substack{-h^{3h-6}(h-1)^{3(n-r)} [\#\mathcal{P}_1^{h}\subseteq \mathcal{T}]^3 \\
        + 3 h^{2h-4} (h-1)^{2(n-r)} [\#\mathcal{P}_1^{h}\subseteq \mathcal{T}]^2 \\
-2h^{h-2}(h-1)^{n-r}[\#\mathcal{P}_1^{h}\subseteq \mathcal{T}]\\   -12h^{2h-4} (h-1)^{n-2r+1}[\#\mathcal{P}_2^{h}\subseteq \mathcal{T}] \\  -6h^{3h-6} (h-1)^{n-3r+2} [\#\mathcal{P}_3^{h}\subseteq \mathcal{T}] \\
 - 12 h^{3h-6}(h-1)^{n-3r+2}   [\#\mathcal{S}_3^{h}\subseteq \mathcal{T}] \\   + 6 h^{3h-6} (h-1)^{2n-3r+1} [\#\mathcal{P}_1^{h}\subseteq \mathcal{T}][\#\mathcal{P}_2^{h}\subseteq \mathcal{T}]}  \right) & \text{ if } i=3h, \\
        0 & \text{ otherwise. }
    \end{cases}
\end{equation*}
\end{theorem}

Fan, Yang and Zheng \cite{fan2024high} have obtained the expression for the order $t$ trace of the adjacency hypermatrix of the power hypergraph of the unicyclic hypergraph in terms of the number of the sub-hypertrees. 

\begin{theorem}\cite{fan2024high}
Let $U:=U_n$ be a unicyclic graph obtained from a cycle $C_{n'}, n' \leq n$ by attaching rooted trees $T_1,\ldots,T_{n'}$, with their roots identified with $n'$ distinct vertices of the cycle. If $U^h$ denotes the $h$-power hypergraph of $U$, then  $Tr_t(U^h)= $
\begin{equation*}
    \begin{cases}
        \substack{ \sum\limits_{k=1}^{n'} n'(h-1)^{(h-1)(n-k)-1}h^{k(h-2)+1}\sum\limits_{T' \in {T(k,2)}} c_{n'}(T')[\#T' \subseteq U] + \\2n'(n'+1)(h-1)^{(h-1)(n-n')}h^{n'(h-2)}} & \text{ if } \frac{t}{h} = n', \\
        \sum\limits_{k=1}^{\frac{t}{h}} t(h-1)^{(h-1)(n-k)-1} h^{k(h-2)} \sum\limits_{T' \in T(k,2)} c_{\frac{t}{h}}(T') [\#T' \subseteq U] & \text{ if } \frac{t}{h} < n'.
    \end{cases}
\end{equation*}
\end{theorem}
In \cite{fan2024trace} Fan, Zheng and Yang have obtained the order $t$ spectral moments (order $t$ traces) of the linear unicyclic hypergraphs.

Girth of a graph $G$, denoted $g(G)$ is the length of the shortest cycle in $G$. 
In \cite{fan2024high}, authors have obtained the expression (in terms of the sub-hypergraphs) for the order $t$ trace of a power hypergraph with some condition on the girth of the graph.

\begin{theorem}\cite{fan2024high}
    Let $G$ be graph on $n'$ vertices and $m$ edges. Then for each $h \geq 3,$ and each $t$ with $h \vert t$ and $\frac{t}{h} < g(G)$,
    \begin{equation*}
        Tr_t(G^h) = \sum\limits_{k=1}^{\frac{t}{h}} t (h-1)^{n'-m-1+(h-1)(m-k)} h^{k(h-2)} \sum\limits_{T' \in T(k,2)} c_{\frac{t}{h}}(T') [\#T' \subseteq G].
    \end{equation*}
\end{theorem}

\subsubsection{Estrada index of Hypergraphs}

If $\lambda_1, \ldots, \lambda_{d}$ are all the eigenvalues of the hypergraph $\mathcal{H}$, where $d = n (h-1)^{n-1}$, then Estrada index of the hypergraph $\mathcal{H}$ denoted by $ EE(\mathcal{H})$ is defined \cite{zhou2023estrada} as $$ EE(\mathcal{H}) = \sum\limits_{j=1}^{d} e^{\lambda_j}. $$
Since $ Tr_t(\mathcal{H}) = \sum\limits_{j=1}^d \lambda_j^t,$ we have $EE(\mathcal{H}) = \sum\limits_{t=0}^\infty  \frac{Tr_t(\mathcal{H})}{t!}. $

The energy of the hypergraph is defined as $E(\mathcal{H})= \sum\limits_{j=1}^d \vert \lambda_i \vert.$
The energy of a simple graph is the highly explored \cite{li2012graph, gutman2001energy,gutman2017survey} area of the spectral graph theory, which was introduced \cite{gutman1978energy} by Gutman (a chemist) in the year 1978. 
There are many articles that consider the matrix associated with the hypergraph and the same invariant is studied \cite{cardoso2020energies,cardoso2022adjacency,kurian2024bounds, fu2024incidence}.

\begin{theorem}\cite{zhou2023estrada}
    Let $\rho$ be the spectral radius of $\mathcal{H}\in  \prescript{n}{m}{H}^{(h)}$, $m \geq 1$. 
    If $\sigma(\mathcal{H})$ is $3$-symmetric, then 
    $$ EE(\mathcal{H}) \leq  \frac{2(\cosh{\rho} -1)}{3 \rho } E(\mathcal{H}) + n2^{n-1} $$
\end{theorem}

\begin{theorem}\cite{zhou2023estrada}
    Let $\mathcal{H}$ be a $h$-uniform hypergraph of order $n$ and size $m$, $m \geq 1.$ Then 
    $$ EE(\mathcal{H}) \geq \frac{m h^{h-2} (h-1)^{n-h-1}}{(h-2)!} + n(h-1)^{n-1}, $$
    with equality if and only if $\mathcal{H}$ is an empty hypergraph.
\end{theorem}

Fan et al. \cite{fan2023trace} have obtained the expression for the order-$t$ trace of the coalescence operation of two hypergraphs in the more general setting. 
By using this result they have characterized the extremal hypertrees that attains maximum and minimum Estrada index among all hypertrees of given size.

\begin{theorem}\cite{fan2023trace}
    Let $\mathcal{T} \in \prescript{}{m}{H}^{(h)}$. 
    Then, 
    $$ EE(\mathcal{P}_m^{(h)}) \leq EE(\mathcal{T}) \leq EE(\mathcal{S}_m^{(h)}), $$
    with equality in the lower bound  is attained if and only if $\mathcal{T} \cong \mathcal{P}_m^{(h)}, $ and the equality in the upper bound is attained if and only if $\mathcal{T} \equiv \mathcal{S}_m^{(h)}. $
\end{theorem}

\begin{theorem}\cite{fan2024trace}
If $\mathcal{T} \in \prescript{}{m}{H}^{(h)}$ has a perfect matching, then 
$$ EE(\mathcal{T}) \leq EE(\mathcal{T}_1),$$
   and the equality is attained if and only if  $\mathcal{T} \cong \mathcal{T}_1$, where the hypertree $\mathcal{T}_1$ is defined in Definition \ref{extremaltree_perfectmatching_spectralradius}.
\end{theorem}

They have also characterized \cite{fan2024trace} the hypertree that attains the maximum Estrada index among the class of all $h$-uniform unicyclic hypergraphs of girth $3$.

\section{Concluding Remarks}

 In this article we have listed most of the pioneering works in the field of the spectral theory of uniform hypergraphs via the normalized adjacency hypermatrix. 
Although many researchers were employing hypermatrices to investigate the spectra of hypergraphs, some were not doing so because there wasn't an effective way to represent general (not uniform) hypergraphs using hypermatrices and also due to the computational complexity in hypermatrices.
 In 2017  Banerjee, Char and Mondal proposed \cite{banerjee2017spectra} a generalization  of the hypermatrix representation of general (not uniform) hypergraphs.  
 Later, many have contributed \cite{banerjee2017spectra,sun2019spectral,ouvrard2017adjacency,zhang2017some,hou2020homogeneous,galuppi2023spectral,zhou2024spectral,duan2021bounds,zhang2022some,shirdel2020non} significantly in obtaining the eigenvalues, their multiplicities, and characterizing the extremal hypergraphs.
 The equitable partition of the graph is a well-studied topic in graph theory, whereas the concepts of equitable partition \cite{taranenko2024perfect} (in general for hypermatrices \cite{jin2018equitable}) and the covering problems \cite{song2024spectral} for studying the spectral property hypergraphs have only been recently initiated.

\section{Backmatter}


\subsection*{Acknowledgment}

Authors thank the organization for the overall support.
Authors also thank Professor Joshua Cooper for his insightful comments and suggestions during the preparation of this manuscript.

\subsection*{Declarations}

\begin{itemize}
\item Funding : Not applicable
\item Conflict of interest: Authors declare that they have no conflict of interest.
\item Ethics approval and consent to participate : Not applicable.
\item Consent for publication : Authors give their consent for publication on acceptance.
\item Data availability : Not applicable.
\item Materials availability : Not applicable.
\item Code availability : Not applicable.
\item Author contribution : Both the authors have equally contributed in writing/reviewing the manuscript.
\end{itemize}

\bibliographystyle{abbrv}
\bibliography{bjourdoc}


\end{document}